\documentclass[11pt,twoside]{article}
\usepackage{latexsym}
\usepackage{amssymb,amsbsy,amsmath,amsfonts,amssymb,amscd}
\usepackage{graphicx,color}
\usepackage{float,url}
\usepackage{cancel}
\usepackage[colorlinks=true]{hyperref}

\newcommand{\commentout}[1]{}

\newcommand{\R}{\mathbb{R}}

\newcommand {\al} {\alpha}

\newcommand {\Chi} {{\bf \raise 2pt \hbox{$\chi$}} }
\newcommand {\f}   {\frac}
\newcommand {\p}   {\partial}
\newcommand{\dis}{\displaystyle}
\newcommand {\proof} {\noindent {\bf Proof}. }
\newcommand{\beq}{\begin{equation}}
\newcommand{\eeq}{\end{equation}}
\newcommand{\bea} {\begin{array}{rl}}
\newcommand{\eea} {\end{array}}
\newcommand{\bepa}{\left\{ \begin{array}{l}}
\newcommand{\eepa} {\end{array}\right.}
\newtheorem{theorem}{Theorem}%[section]
\newtheorem{lemma}[theorem]{Lemma}
\newtheorem{definition}[theorem]{Definition}
\newtheorem{remark}[theorem]{Remark}

\newcommand{\qed}{{ \hfill
                       {\unskip\kern 6pt\penalty 500 \raise -2pt\hbox{\vrule\vbox to 6pt{\hrule width 6pt
                       \vfill\hrule}\vrule} \par}   }}
\title{Reaction-diffusion systems with  initial data of low regularity\\
}

\author{El-Haj LAAMRI\thanks{Universit\'e de Lorraine, Institut Elie Cartan, Email: el-haj.laamri@univ-lorraine.fr}
 \and
Beno\^ \i t PERTHAME\thanks{Sorbonne Universit{\'e}, CNRS, Universit\'{e} de Paris,  Inria, Laboratoire Jacques-Louis Lions UMR7598, F-75005 Paris. 
Email: Benoit.Perthame@sorbonne-universite.fr} 
\thanks{B.P. has received funding from the European Research Council (ERC) under the European Union's Horizon 2020 research and innovation programme (grant agreement No 740623)}
}
\date{\today}

\begin{document}
\maketitle
\pagestyle{plain}
%\tableofcontents
\pagenumbering{arabic}

\begin{abstract}  Models issued from ecology, chemical reactions and several other application fields lead to semi-linear parabolic equations with super-linear growth. Even if, in general, blow-up can occur, these models share the property that mass control is essential. In many circumstances, it is known that this $L^1$ control is enough to prove the global existence of weak solutions. The theory is based on basic estimates initiated by M. Pierre and collaborators, who have introduced methods to prove $L^2$ a priori estimates for the solution.

Here, we establish such a key estimate with initial data in $L^1$ while the usual theory uses $L^2$. This allows us to greatly simplify the proof of some  results. We also establish new existence results of semilinearity which are super-quadratic as they occur in complex chemical reactions. 
Our method can be extended to semi-linear  porous medium equations.
\end{abstract} 
\vskip .7cm

\noindent{\makebox[1in]\hrulefill}\newline

%%%%%%%%%%%%%%%%%%%%%%%%%%%%%%%%%%%%%%%%%%%%
\section*{Introduction}
\label{sec:intro}
%-------------------------------------------
%%%%%%%%%%%%%%%%%%%%%%%%%%%%%%%%%%%%%%%%%%%%

Ecology with Lotka-Volterra systems, chemistry with reaction-rate equations, multi-species diffusion of molecules and many other scientific fields lead to reaction-diffusion systems characterized by different diffusion coefficients. We consider such systems,  set in a smooth domain $\Omega$ of $\R^N$, under the form
\[\bepa
\f{\p}{\p t} u_i(t,x)-  d_i\Delta u_i = f_i(u_1, ...u_m), \qquad x \in \Omega, \; t \geq 0, \; i=1,...,m,
\\[10pt]
\f{\p}{\p \nu} u_i =0 \quad \text{ on } \partial \Omega.
\eepa \]
Even though they are standard, the mathematical understanding of these parabolic systems is still very limited. There are two major obstructions to the construction of weak solutions. Firstly, the right hand sides $f_i$ often has quadratic, and possibly faster, growth for large values of the $u_i$'s. A typical example, for three species is  for $\{i,j,k\} =\{1, 2, 3\}$ (see section~\ref{sec:motivation} for precisions) 
\[
f_i= u_j^\al u_k^\beta - u_i^\gamma
\]
 Secondly, when the diffusion coefficients $d_i$ are very different, there is no maximum principle, and thus a priori estimates are not available besides $L^1$ control.  To circumvent these difficulties, a major $L^1$~theory has been elaborated, which can be roughly summarized as follows: if an $L^1$~bound can be proved for the $f_i$'s, then the existence of a weak solution can be proved, see \cite{Pier2003, Pier2010, LaamPierAHP}. Such an $L^1$~bound can itself be derived thanks to general fundamental lemmas which were initiated by Pierre and Schmitt \cite{PS1997, PS2000}. The most elaborate version is the $L^2$~lemma of Pierre \cite{Pier2010} which we recall below. 

\medbreak
 
 Our main contribution is to construct an existence theory based on the following $L^1\cap H^{-1}$ lemma, which is typically applied to combinations of equations in the above  system. Namely, consider a constant $B\in [0,\infty)$, smooth functions $U, \; F: [0,+\infty)\times \Omega \rightarrow \R^+$ and $V: [0,+\infty)\times \Omega \rightarrow \R^+$such that  $\dis\int_\Omega \Delta_x V(t,x)dx=0$, and satisfying the relations
\begin{equation} 
\tag{R}
\left\{
\begin{array}{ll}
	\f{\p}{\p t} U(t,x)-  \Delta_x V =B -F(t,x),  \quad t\geq 0, \; x \in \Omega \subset \R^N, 
	\\[10pt]
	\f{\p U}{\p \nu} (t,x) =0  \quad \text{on } (0,+\infty)\times\p \Omega,
	\\[10pt]
	U(t=0)=U^0\geq 0  \quad \text{in } \Omega.
\end{array} \right.
\label{eq:mpg1}
\end{equation}
Then, we have
%
%-------------------------------------
\begin{lemma} [First key  estimate with $L^1$ data] 
\label{lm:mpg}
 With the notations above, assume $F\geq 0$, $U\geq 0$, $U^0\in L^1(\Omega)^+\cap H^{-1}(\Omega)$ and $\dis\int_\Omega \Delta_x V(t,x)dx=0$. Then,  there exists a  constant $C>0$ depending only on $\Omega$, such that
\begin{equation}\label{eq:mpg2}
 \int_0^T \int_\Omega  UV   \leq K(T) + \f 12 \| U^0\|_{H^{- 1}(\Omega)}^2 \quad \text{where}\quad   K=  \int_0^T \big[ C \langle F(t) \rangle + \langle V(t) \rangle \big]\; \big[B |\Omega| t+  \int_\Omega U^0 \big] dt .
\end{equation}
 \end{lemma}
 %-------------------------------------
 
The proof of this lemma, and of a more precise version,  is given in Appendix~\ref{ap:KeyProof}, and is reminiscent of the lifting method introduced in~\cite{BloweyElliott91} for the Cahn-Hilliard equation. Several variants can be derived, for instance in dimension $N\leq 3$, a variant allows us to relax the assumption on the sign of~$F$, see Appendix~\ref{sec:estimateA}. Also, Lemma ~\ref{lm:mpg} will be used for the porous media equation with a little more work.
\vskip1mm
This result can be compared with the $L^2$ version of Pierre~\cite{Pier2010} (we indicate here a simple variant).
 %-------------------------------------
\begin{lemma} [Pierre's   $L^2$  lemma] \label{Pierre-lemma}
Assume $V= A U$ with $A(t,x)\geq 0$ and bounded $i.e$ $(0<a\leq A(t,x)\leq b<\infty)$, $U\geq 0$, $U^0\in  L^2(\Omega)$ and \eqref{eq:mpg1} holds, then
\begin{equation}\label{PierreEstimate}
a  \int_0^T \int_\Omega  U^2   \leq  \left\|   \int_0^T A (t,x)dt \right\| _{L^\infty(\Omega)} \int_\Omega  (U^0)^2 .
\end{equation}
 \end{lemma}
%  -------------------------------------
For the ease of the reader, we reproduce the proof in Appendix~\ref{Ap:Pierre}. 
\vskip1mm
This  \textit{a priori $L^2$ estimate} (\ref{PierreEstimate}) was stated in \cite{PS1997, PS2000} and then widely exploited, see for example \cite{DFPV, Pier2010,BLMP, CanDesvFell, PierRoll} and the references therein. Moreover, this estimate has a natural extension to the case where the diffusion operators are degenerate of porous media  type, i.e., $\Delta(u_i^{m_i})$ where $m_i>1$, see \cite[Theorem 2.7]{LaamPierAHP}.
\smallbreak
\noindent Unfortunately, Estimate (\ref{PierreEstimate}) heavily depends  on the $L^2(\Omega)$-norm of the initial data, and does not hold any longer  when the initial data are merely  in $L^1(\Omega)$. However,  the space $L^1(\Omega)$ is a natural functional setting for the systems at hand, as shown by the a priori $L^1(\Omega)$-uniform estimates~(\ref{ControlTotalMass}) and~(\ref{EstimateLUN}) (see section \ref{sec:motivation} below).
\medbreak
The paper is composed of four sections and four appendices.
 In section~\ref{sec:motivation}, we expose in detail the type of nonlinearities $f_i$'s which have usually been treated so far in the literature, and we present the state of the art about existence theory. Section~\ref{sec:Applications1} is devoted to show how Lemma~\ref{lm:mpg} can be used to prove new results, namely, to extend the range of possible nonlinearities $f_i$, or to lower the required regularity on the initial data. Furthermore, our approach simplifies the existing proofs known for $L^1$ initial data. The case of porous media type systems is treated in section~\ref{sec:porous} thanks to an application of  Lemma~\ref{lm:mpg}. Lemma~\ref{lm:mpg} is proved in appendix~\ref{ap:KeyProof} and it is complemented in appendix~\ref{sec:estimateA} with the case where the right hand side $F$ belongs to $L^1$ with no sign condition. Lemma~\ref{Pierre-lemma} is proved in appendix~\ref{Ap:Pierre} while appendix~\ref{Ap:PM} is devoted to the case of the porous media system.
\vskip1mm
We have tried to write the paper in an almost self-contained form, moreover we give precise references for all the points that are not detailed in the work.

%%%%%%%%%%%%%%%%%%%%%%%%%%%%%%%%%%%%%%%%%%%%
%-------------------------------------------
\section{Assumptions, overview and two motivating examples}
\label{sec:motivation}
%-------------------------------------------
%%%%%%%%%%%%%%%%%%%%%%%%%%%%%%%%%%%%%%%%%%%%
Throughout this paper, $\Omega\subset\R^N$ is open bounded with ``good enough" boundary, and we denote
\\
 $Q:=(0,+\infty)\times\Omega, Q_T:=(0,T)\times\Omega$, $\Sigma:=(0,+\infty)\times\partial\Omega$, $\Sigma_T:=(0,T)\times \partial\Omega$ and,  for $p\in [1,+\infty)$
$$\|u(t, .)\|_{L^p(\Omega)}=\left(\int_\Omega |u(t,x)|^p\,dx\right)^{1/p}, \quad \|u\|_{L^p(Q_T)}=\left(\int_0^T\int_\Omega |u(t,x)|^p\,dt dx\right)^{1/p}.$$
%$$\|u(t)\|_{L^\infty(\Omega)}=\text{ess sup}_{x\in \Omega} |u(t,x)|, \quad\|u\|_{L^\infty(Q_T)}=\text{ess sup}_{(t,x)\in Q_T} |u(t,x)|.$$

%--------------------------------------------------------
\subsection{Assumptions} 

 Let us consider a general $m\times m$ reaction-diffusion system
\begin{equation}
\tag{S} \left\{
\begin{array}{lll}
\forall i=1,...,m\\
\frac{\p u_i}{\p t}-d_i\Delta u_i=f_i(u_1, u_2,...,u_m) \;\;\;\;&\text{ in }& \;Q_T,\\[5pt]
\frac{\p u_i}{\p \nu}=0\;\;\;\;&\text{ on }&\;\Sigma_T,\\[5pt]
u_i(0,\cdot)=u_i^0(.)\geq 0 \;\;\;\;& \text{ in } &\;\Omega 
\end{array}
\right.
\label{eq:S}
\end{equation}
%where $\Omega\subset\R^N$ is open bounded with ``good enough" boundary, 
where for all $i=1,\cdots, m$, $d_i>0$ and $f_i:[0,+\infty)^m \to \R$ is locally Lipschitz continuous. 

Moreover, we assume that  the nonlinearities satisfy :\\
\textbf{ (P) :} for all $i=1,\cdots, m$, $ f_i$ is quasi-positive $i.e$ \\
$$\forall \textbf{r}=(r_1,r_2,\cdots,r_m)\in [0,+\infty)^m,\; f_i(r_1\cdot,r_{i-1},0,r_{i+1},...,r_d)\geq 0.
$$
\textbf{ (M) :} there exists $(a_1,\cdots,a_m)\in (0,+\infty)^m$ such that  $\forall \textbf{r}\in [0,+\infty)^m,\; \dis\sum_{1\leq i\leq m}a_if_i(\textbf{r})\leq 0$.
\\
\begin{remark} \label{remarque 1} 
Most of our results can be  extended to:
\\
{\bf (M')} $\forall \textbf{r}\in [0,+\infty[^m,\; \dis\sum_{1\leq i\leq m}a_if_i(\textbf{r})\leq C[1+\dis\sum_{1\leq i\leq m}r_i]$,\\
or to the more general situation  when the  nonlinearities $f_i$ depend also on $t$ and $x$:\\
{\bf (M'')} $\forall (t,x) \in Q$ and $\forall \textbf{r}\in [0,+\infty[^m,\; \dis\sum_{1\leq i\leq m}a_if_i(t,x,\textbf{r})\leq H(t,x) +C\dis\sum_{1\leq i\leq m}r_i$
where $H\in L^1(Q)$.
\end{remark}

%%%%%%%%%%%
\vskip2mm
Properties \textbf{(P)}  and \textbf{(M)} (or \textbf{(M')}) appear naturally in applications. Indeed, evolution reaction-diffusion
systems are mathematical models for evolution phenomena undergoing at the same time
spatial diffusion and (bio-)chemical type of reactions. The unknown functions are generally
densities, concentrations, temperatures so that their nonnegativity is required.
Moreover, often a control of the total mass, sometimes even preservation of the total
mass, is naturally guaranteed by the model. Interest has increased recently for these models
in particular for applications in biology, ecology and population dynamics. We refer to~\cite[Section 2]{Pier2010} for examples of reaction-diffusion systems with properties~\textbf{(P)} and \textbf{(M)} or \textbf{(M')}.
%
%
%---------------------------------
\subsection{Overview}

Many mathematical results are known about the global weak solutions to System~\eqref{eq:S} and we recall some of them and refer to~\cite{QS} for classical solutions. 
\medbreak
\noindent First of all, let us make precise what we mean by solution to~\eqref{eq:S} on $Q_T=(0,T)\times \Omega$.
\vskip1mm
By {\it classical solution} to~\eqref{eq:S}, we mean that, at least\\
$(i)$ $u=(u_1,\cdots,u_m)\in \mathcal{C}([0,T);L^1(\Omega)^m)\cap L^\infty([0,\tau]\times\Omega)^m, \forall \tau \in (0,T)$ ;\\
$(ii)$ $\forall k,\ell=1\dots N$, $\forall p\in (1,+\infty)$
 $$\p_t u_i,\p_{x_k}u_i,\p_{x_kx_\ell}u_i \in L^p(Q_T) \quad i=1, \; \cdots, \; m\,;$$
$(iii)$ equations in~\eqref{eq:S}  are satisfied a.e (almost everywhere).
%%%%%%%%%%%%%%%%%
\vskip2mm
By {\it weak solution} to~\eqref{eq:S} on $Q_T$, we essentially  mean solution in the sense of distributions or, equivalently here, solution in the sense of Duhamel's formula with the corresponding semigroups. More precisely,  for $1\leq i\leq m$, $f_i(u)\in L^1(Q_T)$ and 
\begin{eqnarray*}
u_i(t, .)&=&S_{d_i}(t)u_i^0(.) + \int_0^t S_{d_1}(t-s)f_i(u(s, .))\,ds
\end{eqnarray*}
where $S_{d_i}(.)$ is the semigroup generated in $L^1(\Omega)$  by $-d_i\Delta$  with homogeneous Neumann boundary condition.
\vskip2mm
%%%%%%%%%%%%%%%%%%%%%%%%%%%%%%%%%%%%%%%%%%%%%%%%%%%%%%%%%%%%%%%%ù
%
For initial data $u^0\in (L^\infty(\Omega))^m$, the local Lipschitz continuity  of the nonlinearities implies the existence of a local classical solution to~\eqref{eq:S} on a maximal interval $[0, T_{\max})$. Moreover,  the initial data are nonnegative,
so the quasi-positivity $\textbf{(P)}$ ensures that the solution stays nonnegative as long as it exists, see~\cite{QS}.
The assumption $\textbf{(M)}$ gives an upper bound on the total mass of the system $i.e$ for all $t \in (0, T_{\max})$
\begin{equation} \label{ControlTotalMass}
\sum_{i=1}^m \int_\Omega a_ iu_i(t,x)\,dx \leq  \sum_{i=1}^m\int_\Omega a_iu_i^0(x)\,dx.
\end{equation}
In fact, multiplying each $i$-th equation by $a_i$ and adding the $m$
equations to obtain
\begin{equation}\label{Somme}
\sum_{i=1}^m a_i\p_tu_i-\sum_{i=1}^m a_id_i\Delta u_i = \sum_{i=1}^m a_if_i\leq 0.
\end{equation}
By integrating (\ref{Somme})  on $(0,t)\times\Omega$ and taking into account the boundary conditions $\dis\int_\Omega \Delta u_i=\dis\int_{\p\Omega} \frac{\p u_i}{\p\nu}=0$ and  $\textbf{(M)}$, we obtain (\ref{ControlTotalMass}).
In other words, the total mass of $m$  components is preserved.
Together with the nonnegativity of $u_i$, estimate (\ref{ControlTotalMass}) implies that
\begin{eqnarray}\label{EstimateLUN}
\forall t\in (0,T_{\max})\;,\; \|u_i(t,.)\|_{L^1(\Omega)}\leq \frac{1}{a_i}\|\sum_{j=1}^m a_ju_j^0\|_{L^1(\Omega)}.
\end{eqnarray}
So the $u_i(t,.)$ remain bounded in $L^1(\Omega)$  uniformly in time as long as solution exists.
\vskip1mm
Let us emphasize that  uniform $L^\infty$-bounds, rather than $L^1$-bounds, would provide global
existence in time of \textit{classical} solutions, by the standard theory for reaction-diffusion systems. 
The point here is that bounds are a priori only in $L^1$ and one cannot apply the 
$L^\infty$-approach (see~\cite{QS}) even if the initial data are regular except in the restrictive  case where $d_1=\cdots=d_m$.\\ However the situation is quite more
complicated if the diffusion coefficients are different from each other. Some additional assumptions on the structure of the source terms are needed for global existence of classical  solutions.
This question has   been  widely  studied.
For some recent results, see e.g \cite{LaamAAM, FellLaam, GV10, PierSuzYam2019, SuzYam, CapGoudVass,  Soup18, FellMorgTang}. Concerning the results established before 2010,  we refer the intrested reader to the exhaustive survey \cite{Pier2010} for a general presentation of the problem, further references, and many deep comments on the mathematical difficulties raised by such systems. In fact,  the solutions may blow up
in $L^\infty(\Omega)$ in finite time as proved in \cite{PS1997, PS2000}  where explicit finite time $T^*$ blow up is given. 
More precisely, Pierre and Schmitt exhibited a system with two species fulfilling $\textbf{(P)}$ and $\textbf{(M)}$ with $d_1\neq d_2$ and  strictly superquadratic polynomial nonlinearities $f_i$ such that  $T^*<+\infty$ and
\begin{equation*}
\lim_{t\nearrow T^*}\Vert u_1(t, .)\Vert_{L^\infty(\Omega)}= \lim_{t\nearrow T^*}\Vert u_2(t, .)\Vert_{L^\infty(\Omega)}=+\infty.
\end{equation*}

Thus, even in the semi-linear case  it is necessary to deal  with weak solutions if one expects global existence for more general nonlinearities and initial data of low regularity. 
\vskip 1mm
It is worth to pointing out that, while considerable effort has been devoted to the study of systems with initial data in $L^\infty(\Omega)$ or~$L^2(\Omega)$ and at most quadratic nonlinearities, relatively little is known in the case of systems with initial data in~$L^1(\Omega)$. We refer to \cite{{Pier2003,Pier2010}} for linear diffusion, and \cite{Laamthese, LaamPierAHP, LaamPierM3AS, PierRoll} for nonlinear diffusion.
\vskip1mm
\noindent Here, we just  recall the  following three theorems closely related to our present study.
\vskip1mm
\noindent $\bullet$ \textbf{$L^1$-Theorem:}\\
In order to give the precise statement, let us  introduce  the following approximation of System~\eqref{eq:S}
\begin{equation}\label{Sn}
\left\lbrace\begin{array}{lll}
i=1,...,m,\\[5pt]
\dis\partial_tu_i^n-\Delta u_i^n=f_i^n(u^n) &\text{ in } & Q=(0,+\infty)\times\Omega, \\[5pt]
\frac{\p u_i^n}{\p \nu}(t,.) =  0 &\text{ on }& \Sigma=(0,+\infty)\times\p \Omega,\\[5pt]
u_i^n(0,.)=u_{i,0}^n\geq 0 & \text{ in } & \Omega ,\\
\end{array}
\right.
\end{equation}
where $u_{i,0}^n:= \inf\{u_{i,0},\; n\}$ and  $f_i^n:=\dis\frac{f_i}{1+\frac{1}{n}\sum_{1\leq j\leq m}|f_j|}$.\\
 For $i=1,\cdots,\,m$, $u_{i,0}^n\in L^\infty(\Omega)$ and converges to $u_{i,0}$ in $L^1(\Omega)$, $f_i^n$ is locally Lipschitz continuous and satisfy \textbf{(P)} and \textbf{(M)}. Moreover
$\|f_i^n\|_{L^\infty(\Omega)}\leq n$. 
Therefore,  the approximate system (\ref{Sn}) has a nonnegative classical global solution $u^n= (u_1^n,\cdots,u_m^n)$ (see e.g  \cite{LSU}).
\begin{theorem} [Pierre, \cite{Pier2010}] \label{Lun} 
Besides {\bf (P)+(M)}, assume that the following \textbf{ a priori $L^1$-estimate} holds: there exists $C(T)$ independent of $n$ such that
\begin{equation} \label{Lun-estimate}
\forall i=1,...,m,\; \forall T>0,\; \int_{Q_T}|f_i^n(u_1^n,\cdots,u_m^n)|\leq C(T).
\end{equation} 
Then, as $n \to \infty$, up to a subsequence,  $u^n$ converges in $L^1(Q_T)$ for all $T>0$ to some {\em global weak solution} $u$ of~\eqref{eq:S} for all $(u_{0,1},\cdots,u_{0,m})\in (L^1(\Omega)^+)^m$. 
\end{theorem}
For a proof, see \cite[Theorem 5.9]{Pier2010}.
%------------------------------------------
\vskip2mm
The two following theorems  are  more adequate for cases where the nonlinearities are at most quadratic or super-quadratic respectively. 
\vskip1mm
%-------------------------------------------------------------
\noindent $\bullet$ \textbf{$L^2$-Theorem with linear diffusion:}
\begin{theorem}[Pierre, \cite{Pier2010}]\label{LdeuxL}
Besides \textbf{(P)} and \textbf{(M)}, assume that the  $f_i$ are at most quadratic $i.e$  theres exists $C>0$ such that for all $i=1,\cdots,m$:

\textbf{\em (QG)} $\quad\quad |f_i(r)|\leq C\big(1+\dis\sum_{j=1}^m r_j^2\big).$ \\

\noindent Then, there exists a {\em global weak solution} to~\eqref{eq:S}
for all  $u_0=(u_{0,1},\cdots,u_{0,m}) \in (L^2(\Omega)^+)^m$. 
\end{theorem} 
For a proof, see \cite[Theorem 5.11]{Pier2010}.
\vskip1mm
%-------------------------------------------------------------
\noindent $\bullet$ \textbf{$L^2$-Theorem with nonlinear diffusion:}
Consider the following $m\times m$ reaction-diffusion system
\begin{eqnarray*}
(NLDS) \left\{
\begin{array}{lll}
\forall i=1,...,m\\
\frac{\p u_i}{\p t}-d_i\Delta (u_i^{m_i})=f_i(u_1, u_2,...,u_m) \;\;\;\;&\text{ in }& \;Q_T,\\[5pt]
\frac{\p u_i}{\p \nu}=0\;\;\;\;&\text{ on }&\;\Sigma_T,\\[5pt]
u_i(0,\cdot)=u_i^0(.)\geq 0 \;\;\;\;& \text{ in } &\;\Omega 
\end{array}
\right.
\end{eqnarray*} 
where for all $i=1,\cdots, m$, $d_i>0$, $m_i>1$  and $f_i:[0,+\infty)^m \to \R$ is locally Lipschitz continuous. 

\begin{theorem}[Laamri--Pierre, \cite{LaamPierAHP}]\label{LdeuxNLN}
Besides \textbf{(P)} and \textbf{(M)}, assume that the  $f_i$'s are at most super-quadratic $i.e$  theres exists $C>0$ such that for all $i=1,\cdots,m$:

\textbf{\em (SQG)} $\quad\quad |f_i(r)|\leq C\big(1+\dis\sum_{j=1}^m r_j^{m_i+1-\varepsilon}\big).$ \\
\noindent Then, there exists a {\em global weak solution} to $(NLDS)$
for all  $u_0=(u_{0,1},\cdots,u_{0,m}) \in (L^2(\Omega)^+)^m$. 
\end{theorem} 
For a proof, see \cite[Theorem 2.7]{LaamPierAHP}.

%%%%%%%%%%%%%%%%%%%
\subsection{Motivating examples}

 Now, let us introduce two systems  we are considering in this paper. In fact, these systems contain
 the major difficulties encountered in a large class of similar problems as regards global existence in time of solutions.
 \\
$\bullet$ \textbf{First system:} Let $(\alpha,\beta,\gamma  )\in [1,+\infty)^3$  and  the  following reaction-diffusion system
$$
(S_{\alpha\beta\gamma}) \quad \left\lbrace\begin{array}{lllll}
\dis\f{\p u_1}{\p t} (t,x)-  d_1\Delta u_1 &=& \alpha\big(u_3^\gamma - u_1^\al u_2^\beta\big) &\text{ in } & Q_T,
	\\[7pt]
	\dis\f{\p u_2}{\p t} (t,x)-  d_2\Delta u_2 &=& \beta\big(u_3^\gamma - u_1^\al u_2^\beta\big) &\text{ in } & Q_T,
	\\[7pt]
	\dis\f{\p u_3}{\p t} (t,x)-  d_3\Delta u_3 &=&\gamma\big( -u_3^\gamma +u_1^\al u_2^\beta\big) &\text{ in } & Q_T,
	\\[5pt]
	\f{\p u_i}{\p \nu} (t,x) &=&0  &\text{ on } & \Sigma_T, 
	%\quad i=1, \; 2, \; 3,
	\\[7pt]
        u_i(t=0)&=&u_i^0 \geq 0 &\text{ in } & \Omega, 
        \\
         i=1, \; 2, \; 3.	
\end{array}
\right.$$
If $\alpha,\beta$ and  $\gamma$ are positive integers, system $(S_{\alpha\beta\gamma})$ is intended to describe  for example the evolution of a reversible chemical reaction of type
 $$\alpha U_1+\beta U_2\rightleftharpoons \gamma U_3$$
where $u_1$, $u_2$, $u_3$ stand for the density of $U_1$,   $U_2$ and $U_3$ respectively.\\
 Let us recall the known results about the global existence of solutions  to system $(S_{\alpha\beta\gamma})$ ; for more details see \cite{LaamAAM} and the references therein.\\
  %It obvious when  $d_1=d2=d_3$. However, it is far of being 
%Indeed,  $Z=u+v+2w$ satisfies
%$$(E)\left\lbrace\begin{array}{llll}
%Z_t-d\Delta Z & = & 0& (0,+\infty)\times\Omega,\\
%%\Part{Z}{n}   & = & 0 & (0,+\infty)\times\partial\Omega,\\
%Z(0,x) & = & Z_0(x) & x\in\Omega,\\
%\end{array}
%\right.$$
%where $Z_0(x)=u_0(x)+v_0(x)+2w_0(x)$.\\
%In particular, we deduce by maximum principle that
%$$\|u(t)+v(t)+2w(t)\|_{L^\infty(\Omega)}\leq \|u_0+v_0+2w_0\|_{L^\infty(\Omega)},\quad\quad t\geq 0.$$
%Together with nonnegativity, this implies
% $$\|u(t)\|_{L^\infty(\Omega)}+ \|v(t)\|_{L^\infty(\Omega)}+ \|w(t)\|_{L^\infty(\Omega)}\leq \|u_0+v_0+2w_0\|_{L^\infty(\Omega)}, \quad\quad t\geq 0.$$
%In other words, $u(t)$, $v(t)$ and $w(t)$ stay uniformly bounded in $L^\infty(\Omega)$ and therefore $T_{\max}=+\infty$.
%\vskip2mm
In the case where the diffusion coefficients are different from each other, global existence   is more complicated (it obviously holds if $d_1=d_2=d_3$). It  has been studied by several authors in the following cases.\\
\textbf{\em First case } $\alpha=\beta=\gamma=1$. In this case, global existence of classical solutions has been obtained by Rothe \cite{Rot} for dimension $N\leq 5$. Later, it has first been proved by Martin-Pierre \cite{MP1} for all dimensions $N$ and then by Morgan  \cite{Mo89}.
\\
 \noindent Global existence of weak solutions  has been proved by Laamri \cite{Laamthese} for initial data only in $L^1(\Omega)$.\\
\noindent\textbf{\em Second case } $\gamma = 1$ regardless of $\alpha$ and $\beta$. In this case, global existence of classical solutions has been obtained by Feng~\cite{F} in  all dimensions $N$ and  more general boundary conditions.\\
\noindent\textbf{\em Third case } $\alpha+\beta< \gamma$, or when $1<\gamma<\dis\frac{N+6}{N+2}$ regardless of $\alpha$ and $\beta$. In these cases, global existence of classical solutions was established  by the first author in \cite{LaamAAM}.\\
Up to our knowledge, in the case $\alpha+\beta \geq \gamma>1$, global existence of classical solutions  remain an open question  when the diffusivities $d_i$ are away from each other.\\
Exponential decay towards equilibrium has been studied by Fellner-Laamri \cite{FellLaam}. \\
\noindent\textbf{\em Fourth case } $\alpha=\beta= 1$ or $\gamma\leq 2$. In this  case, Pierre \cite{Pier2010} has proved global existence of weak solutions for initial data in $L^2(\Omega)$.
\vskip2mm
 
% In this paper, we are able to prove global existence of weak solutions in the following cases
%  ($\gamma \leq 2$ and whatever are $\alpha, \, \beta$) and ($\alpha= \beta=1$ and whatever is $\gamma$) for initial data only in $L^1(\Omega)\cap H^{-1}(\Omega) $.
%  See subsection \ref{sec:alphabetagamma}
%  %
%\medbreak
%
\noindent $\bullet$ \textbf{Second system:} Let us consider 
the following system  naturally arising in chemical kinetics when modeling the following reversible reaction
$$\alpha U_1+\beta U_3\rightleftharpoons \gamma U_2+ \delta U_4$$
$$( S_{\alpha\beta\gamma\delta})\left\lbrace\begin{array}{lllll}
\dis\f{\p u_1}{\p t} (t,x)-  d_1\Delta u_1 &=& \alpha\big(u_2^\gamma u_4^\delta - u_1^\al u_3^\beta\big) &\text{ in } & Q_T,
	\\[5pt]
	\dis\f{\p u_2}{\p t} (t,x)-  d_2\Delta u_2 &=& \beta\big(u_1^\al u_3^\beta-u_2^\gamma u_4^\delta \big) &\text{ in } & Q_T,
	\\[5pt]
	\dis\f{\p u_3}{\p t} (t,x)-  d_3\Delta u_3 &=&\gamma\big( u_2^\gamma u_4^\delta - u_1^\al u_3^\beta\big) &\text{ in } & Q_T,
	\\[5pt]
	\dis\f{\p u_4}{\p t} (t,x)-  d_4\Delta u_4 &=& \delta\big(u_1^\al u_3^\beta-u_2^\gamma u_4^\delta \big) &\text{ in } & Q_T,
	\\[5pt]
	\f{\p u_i}{\p \nu} (t,x) &=&0  &\text{ on } & \Sigma_T, 
	%\quad i=1, \; 2, \; 3,
	\\[5pt]
        u_i(t=0)&=&u_i^0 \geq 0 &\text{ in } & \Omega, 
        \\
         i=1, \; 2, \; 3, \; 4.	
\end{array}
\right.$$

Up to our best knowledge, global existence in the  \textbf{more restrictive case  $\alpha=\beta=\gamma=\delta=1$}  has been  studied by many authors.  Concerning global existence of  weak solution in all dimensions $N$, it  has been  first obtained  in \cite{DFPV} for initial data $u_i^0$ such $u_i^0\log u_i^0\in L^2(\Omega)$, then  in  \cite{Pier2010} for initial data in $L^2(\Omega)$ and later in  \cite{PierRoll} for initial data in $L^1(\Omega)$; see also \cite{PierSuzYam2019, SuzYam} and the references therein.
% About global classical solution, in small dimensions $N\leq2$, it has been first established in \cite{GV10}, then in \cite{CanDesvFell}.  then  {CapGoudVass,  Soup18, FellMorgTang}
In small dimensions $N\leq2$, global existence of classical solutions  is established in \cite{GV10, CanDesvFell, PierSuzYam2019}. Recent results in  \cite{CapGoudVass,  Soup18, FellMorgTang} have proved that global  classical solutions exist in all dimensions.
%
%--------------------------------------
\subsection{Brief Summary}
%----------------------------------------
In this work,  we exploit the ``good" $L^1$-framework provided by the two conditions $\textbf{(P)}$ and $\textbf{(M)}$ on the one hand  and our key   estimate with $L^1$ data (\ref{eq:mpg2}) on the other hand. This allows us to extend known global existence results for initial data in $L^2(\Omega)$ to cases where the initial data only belong to  $L^1(\Omega)\cap H^{-1}(\Omega)$.\\
More precisely, the main results of this paper can be summarized in the following points.\\
1)  We are able to prove global existence of weak solutions
\begin{itemize}
\item  to system $(S_{\alpha\beta\gamma})$ in the following cases
  ($\gamma \leq 2$ and whatever are $\alpha, \, \beta$) and ($\alpha= \beta=1$ and whatever is $\gamma$). 
  %for initial data only in $L^1(\Omega)\cap H^{-1}(\Omega) $.
  See subsection~\ref{sec:alphabetagamma} ;
\item to system $(S_{\alpha\beta\gamma\delta})$  in the following cases
  ($\gamma=\delta=1$ and whatever are $\alpha, \, \beta$) and ($\alpha= \beta=1$ and whatever are $\gamma, \, \delta$).
  % for initial data only in $L^1(\Omega)\cap H^{-1}(\Omega)$.
  See subsection~\ref{sec:alphabetagammadelta}.
  \end{itemize}
2) We  are also able to establish that Theorem~\ref{LdeuxL}  and Theorem~\ref{LdeuxNLN} hold for initial data  only in  $L^1(\Omega)\cap H^{-1}(\Omega)$. See subsection~\ref{sec:quadratic1} and section~ref{sec:porous} respectively.
\vskip2mm
As we will see in more detail in the next sections, our demonstrations are based on two ingredients~: our key   estimate with $L^1$ data (\ref{eq:mpg2}) and  $L^1$-Pierre's theorem, Theorem~\ref{Lun}.
\vskip2mm
  Before ending this section, note that we state our theorems for Neumann boundary conditions, but they can easily be adapted to Dirichlet  boundary conditions. One must however be careful when
choosing two different boundary conditions for $u_i$ and $u_j$ for $i\neq j$, see \cite{BL, MP2}.

%%%%%%%%%%%%%%%%%%%%%%%%%%%%%
%--------------------------------------------------------------------------
\section{Applications of the key $L^1$ estimates}
\label{sec:Applications1}
%--------------------------------------------------------------------------
%%%%%%%%%%%%%%%%%%%%%%%%%%%%%

Our first estimate (\ref{eq:mpg2}) can be exploited for establishing global existence -in- time of weak solutions to a large class of  reaction-diffusion systems  \textit{with initial data of low regularity}.
%%%%%%%%%%%%%%%%%%%%%%%%%%%%%%%%%%%%%%%%%%%%%
\subsection{A classical quadratic model with mass dissipation }
\label{sec:quadratic1}
%-------------------------------------------
%%%%%%%%%%%%%%%%%%%%%%%%%%%%%%%%%%%%%%%%%%%%
As a first use of our estimate (\ref{eq:mpg2}), see also Appendix~\ref{ap:KeyProof}, 
we show here how one may prove global existence of weak solutions for quadratic nonlinearities  and initial data of low regularity. In fact, this case  is of interest due to its relevance in many applications such as chemical reactions or parabolic Lotka-Volterra type systems, see~\cite{BPerth}.
\\
We consider the following general system: for $t\geq 0, \; x \in \Omega \subset \R^N$, $1\leq i\leq m$,
\begin{equation} \left\{
\begin{array}{ll}
	\dis\f{\p}{\p t} u_i(t,x)-  d_i\Delta u_i = f_i(u_1,\cdots,u_m), 
	\\[5pt]
	\dis\f{\p u_i}{\p \nu} (t,x) =0  \text{ on } (0,+\infty)\times\p \Omega,
	\\[5pt]
	u_i(t=0)=u_i^0 \geq 0 ,
\end{array} \right.
\label{eq:quad}
\end{equation}
where  $d_i>0$ and $f_i:\R^m\rightarrow \R$ is locally Lipschitz continuous. Besides  \textbf{(P)} and \textbf{(M)}, assume  that the $f_i$ are  at most quadratic $i.e$ there exists $C>0$ such that for all $ 1\leq i\leq m$\\

%-------------------------------------------------------
\begin{theorem}\label{quadgrow} %%
Besides  \textbf{(P)} and \textbf{(M)}, assume  that the $f_i$ are  at most quadratic $i.e$ there exists $C>0$ such that for all $ 1\leq i\leq m$\\
  \textbf{\em (QG)} $\quad\quad |f_i(r)|\leq C\big(1+\dis\sum_{j=1}^m r_j^2\big).$ \\
 Then,  for all $u^0=(u^0_1,\cdots, u^0_m)\in (L^1(\Omega)^+\cap H^{- 1}(\Omega))^m$, System (\ref{eq:quad}) has a non-negative  weak solution which satisfies for all $T>0$
$$
u_i \in L^2(Q_T), \qquad u_i \in L^\infty \big((0,T); H^{- 1}(\Omega) \big).
$$
\end{theorem} 
\proof 
We approximate  the initial data and right hand side of System (\ref{eq:quad})  with  $u_{i,0}^n:= \inf\{u_{i,0},\; n\}$ and  $f_i^n:=\dis\frac{f_i}{1+\frac{1}{n}\sum_{1\leq j\leq m}|f_j|}$ and set
\begin{equation}\label{Snquad}
\left\lbrace\begin{array}{lll}
i=1,...,m,\\
\f{\p}{\p t}u_i^n-\Delta u_i^n=f_i^n(u^n) &\text{ in } & Q=(0,+\infty)\times\Omega, 
\\[5pt]
\frac{\p u_i^n}{\p \nu}(t,.) =  0 &\text{ on }& \Sigma=(0,+\infty)\times\p \Omega,
\\[5pt]
u_i^n(0,.)=u_{i,0}^n\geq 0 & \text{ in } & \Omega.
\end{array}
\right.
\end{equation}
 For $i=1,\cdots,\,m$, $u_{i,0}^n\in L^\infty(\Omega)$ and converges to $u_{i,0}$ in $L^1(\Omega)$, $f_i^n$ is locally Lipschitz continuous and satisfies  $(P)$ and $(M)$. Moreover
$\|f_i^n\|_{L^\infty(\Omega)}\leq n$. 
Therefore,  the approximate system (\ref{Snquad}) has a nonnegative classical  global solution  (see e.g  \cite{LSU}).\vskip2mm
\noindent Thanks to Theorem~\ref{Lun},  it is sufficient to establish the following a priori $L^1$-estimate
\begin{equation} \label{Lun-bound1}
\forall i=1,...,m,\; \forall T>0,\; \int_{Q_T}|f_i^n(u_1^n,\cdots,u_m^n)|\leq C(T)
\end{equation} 
where $C(T)$ is independent of $n$.\\
Multiplying each $i$-th equation of System (\ref{Snquad}) by $a_i$ and adding the $m$
equations to obtain
\begin{equation}\label{somme2}
\f{\p}{\p t} [\sum_{i=1}^m a_iu_i^n]-\sum_{i=1}^m a_id_i\Delta u_i^n = \sum_{i=1}^m a_if_i^n\leq 0
\end{equation}
Applying Lemma \ref{lm:mpg} with $U=\dis\sum_{i=1}^ma_iu_i^n$, $V=\dis\sum_{i=1}^md_ia_iu_i^n$, $B=0$ and $F=- \dis\sum_{i=1}^m a_if_i^n$. Then,  there exists $C(T)$ independent of $n$ such that
$$\int_0^T \int_\Omega [\sum_{i=1}^ma_iu_i^n] \; [\sum_{i=1}^md_ia_iu_i^n]dtdx\leq C(T).$$
Thanks to their non-negativity,  the $u_i^n$'s  are bounded  in $L^2(Q_T)$ and then the $f_i^n(u^n)$'s are bounded in $L^1(Q_T)$ independently of $n$.
\qed
\begin{remark}
The above proof seems very simple. But it rests on the very powerful Theorem~\ref{Lun}, whose proof is delicate. Moreover, Lemma~\ref{lm:mpg}  directly provides the $L^1$-bound of the nonlinearities~$f_i$.
\end{remark}
\begin{remark} 
With a very fine analysis, Theorem \ref{quadgrow} was proved by Pierre and Rolland in \cite{PierRoll} under an \textit{\em extra assumption} on the growth of the negative part of the reaction terms, namely:\\
$(H)$\; $\exists \varphi \in \mathcal{C}([0,+\infty), [0,+\infty))$ such that for all $r=(r_1,\cdots, r_m)\in [0,+\infty)^m$ and for all $1\leq i\leq m$
$$[f_i(r)-f_i(\pi_i(r))]^-\leq \varphi(r_i)\big[1+\sum_{i=1}^mr_i\big]$$
where $\pi_i(r):= (r_1,\cdots,r_{i-1},0, r_{i+1},\cdots,r_m)$.\\
Also, it was noticed in \cite{PierRoll} that  the nonlinearities $f_2(r_1,r_2)=-f_1(r_1,r_2)= (r_1^2+ r_2^2)\sin(r_1r_2)$ do not satisfy $(H)$: for this example, global existence with general $L^1(\Omega)$ initial data had remained so far an open question. But our theorem does apply to this example.
\end{remark} %
\vskip 1mm
%%-------------------------------------------
%%%%%%%%%%%%%%%%%%%%%%%%%%%%%%%%%%%%%%%%%%%%%
\textbf{Example:} Let us mention that our theorem \ref{quadgrow} applies to the famous Lotka-Volterra's system where 
\begin{equation}
f_i(u)=\left(e_i+\dis\sum_{j=1}^ma_{ij}u_j\right)u_i, \quad 1\leq i\leq m
\end{equation}
in~\eqref{eq:S} and satisfying, with $A:=(a_{ij})\in \mathcal{M}_m(\R)$, 
$$
(e_1,\cdots,e_m)\in ]-\infty,0]^m \text { and } \langle Au,u\rangle \leq 0 \quad\quad \forall (u_1,\cdots,u_m)\in [0,+\infty)^m .
$$
Let us recall that this result was obtained in \cite{SuzYam} under the strong condition
$$^tA=-A$$
and initial data in $(L^2(\Omega)^+)^m$.

%%%%%%%%%%%%%%%%%%%%%%%%%%%%%%%%%%%%%%%%%%%%
\subsection{ A $3\times 3$ system with nonlinearities of general growth}
\label{sec:alphabetagamma}
%-------------------------------------------
%%%%%%%%%%%%%%%%%%%%%%%%%%%%%%%%%%%%%%%%%%%%
To show how far our approach may be carried out, we now consider systems with  nonlinearities of higher degree of the following form:
 for $t\geq 0, \; x \in \Omega \subset \R^N, $
\begin{equation} (S_{\alpha\beta\gamma}) \quad
\left\{
\begin{array}{ll}
	\f{\p}{\p t} u_1(t,x)-  d_1\Delta u_1 = 
	\alpha(u_3^\gamma - u_1^\alpha u_2^\beta) ,
	\\[5pt]
	\f{\p}{\p t} u_2(t,x)-  d_2 \Delta u_2 = \beta(u_3^\gamma - u_1^\alpha u_2^\beta) ,
	\\[5pt]
	\f{\p}{\p t} u_3(t,x)- d_3 \Delta u_3 =  \gamma(u_1^\alpha u_2^\beta -u_3^\gamma) ,
	\\[5pt]
	\f{\p u_i}{\p \nu} (t,x) =0  \text{ on } (0,+\infty)\times\p \Omega,
	\\[5pt]
        u_i(t=0)=u_i^0 \geq 0 ,
        \\[5pt]
        1\leq i\leq 3
\end{array} \right.
\label{eq:alphabetagamma}
\end{equation}
where $(\alpha, \beta, \gamma)\in [1,+\infty[^3$.
\begin{theorem} 
Assume $u^0_i \in L^1(\Omega)^+\cap H^{- 1}(\Omega)$, $i=1, \; 2, \; 3$.\\
1) 
For $\gamma\leq 2$ and any $\al, \; \beta$, System (\ref{eq:alphabetagamma}) has a non-negative weak solution which 
satisfies, for all $T>0$,
\[ u_i \in L^2(Q_T), \quad u_i \in L^\infty \big((0,T); H^{- 1}(\Omega) \big), \quad
\int_0^T \int_\Omega u_1^\alpha u_2^\beta dx dt \leq C(T).
\]
2) 
For $\al= \beta =1$ and any $\gamma$, System (\ref{eq:alphabetagamma}) has a non-negative weak solution which 
satisfies, for all $T>0$,
\[ u_i \in L^2(Q_T), \quad u_i \in L^\infty \big((0,T); H^{- 1}(\Omega) \big), \quad
\int_0^T \int_\Omega u_3^\gamma dx dt \leq C(T).
\]
\end{theorem} 

\proof
For the  simplicity of the presentation, we consider the truncated system as before, but we drop the indexation by $n$.
 Let $T>0$. Using the conservation laws for $\gamma u_1 +\alpha u_3$ and $\gamma u_2+\beta u_3$ and Lemma~\ref{lm:mpg}, 
we conclude that $u_i$ is bounded in $L^2(Q_T)$.\\
\textit{First case} $\gamma \leq 2$.
We integrate the first (or the second) equation to obtain
\[
\int_\Omega u_1(T,x)  dx  + \alpha\int_0^T \int_\Omega u_1^\alpha u_2^\beta dx dt = \alpha\int_0^T \int_\Omega u_3^\gamma  dx dt  + \int_\Omega u_1^0(x)  dx.
\]
All the terms on the right-hand side are bounded, we conclude the  $L^1$-bound on $ u_1^\alpha u_2^\beta $.\\
\textit{Second case} $\alpha =\beta= 1$.
We integrate the third equation to obtain
\[
\int_\Omega u_3(T,x)  dx  + \gamma\int_0^T \int_\Omega u_3^\gamma  dx dt =  \int_\Omega u_3^0(x)  dx + \gamma\int_0^T \int_\Omega u_1 u_2dx dt.
\]
All the terms on the right hand side being bounded, we conclude the  $L^1$-bound on $u_3^\gamma $.
\vskip2mm
\noindent In the two cases, $u_1^\alpha u_2^\beta -u_3^\gamma$ is  bounded  in $L^1(Q_T)$. Therefore, we conclude thanks to Theorem~\ref{Lun}.\qed

\begin{remark}  
\label{rmk-alphabetagamma}
When $\gamma\leq 2$,  the specific form $u_1^\alpha u_2^\beta$ does not play any role, and any function $f(u_1,u_2) \geq 0$ satisfying $f(0,u_2)= f(u_1,0)=0$ will work.\\
In the same way, when $\alpha =\beta= 1$,  the specific form $u_3^\gamma$ does not play any role, and any function $g(u_3) \geq 0$ satisfying $g(0)=0$ will work.
\end{remark}
%%%%%%%%%%%%%%%%%%%%%%%%%%%%%%%%%%%%%%%%%%%%
\subsection{ A $4\times 4$  system with nonlinearities of general growth}
\label{sec:alphabetagammadelta}
%-------------------------------------------
%%%%%%%%%%%%%%%%%%%%%%%%%%%%%%%%%%%%%%%%%%%%
In this subsection, let us consider the following system  naturally arising in chemical kinetics when modeling the following reversible reaction
$$\alpha U_1+\beta U_2\rightleftharpoons \gamma U_3+ \delta U_4.$$
 For $t\geq 0, \; x \in \Omega \subset \R^N, $
\begin{equation} (S_{\alpha\beta\gamma\delta}) \quad 
\left\{
\begin{array}{ll}
	\f{\p}{\p t} u_1(t,x)-  d_1\Delta u_1 = \alpha( u_3^\gamma u_4^\delta - u_1^\alpha u_2^\beta) ,
	\\[5pt]
	\f{\p}{\p t} u_2(t,x)-  d_2 \Delta u_2 = \beta( u_3^\gamma u_4^\delta - u_1^\alpha u_2^\beta) ,
	\\[5pt]
	\f{\p}{\p t} u_3(t,x)- d_3 \Delta u_3 =  \gamma(u_1^\alpha u_2^\beta -u_3^\gamma u_4^\delta) ,
	\\[5pt]
	\f{\p}{\p t} u_4(t,x)- d_4 \Delta u_4 =  \delta(u_1^\alpha u_2^\beta -u_3^\gamma u_4^\delta) ,
	\\[5pt]
	\f{\p u_i}{\p \nu} (t,x) =0  \text{ on } (0,+\infty)\times\p \Omega,
	\\[5pt]
        u_i(t=0)=u_i^0 \geq 0 ,
        \\[5pt]
        1\leq i\leq 4.
\end{array} \right.
\label{eq:alphabetagammadelta}
\end{equation}

\begin{theorem} \label{thm:alphabetagammadelta}
Assume $u^0_i \in L^1(\Omega)^+\cap H^{- 1}(\Omega)$, $i=1, \; 2, \; 3,\; 4$.\\
1) 
For $\gamma=\delta = 1$ and whatever are $\al, \; \beta$, System (\ref{eq:alphabetagammadelta}) has a non-negative weak solution which 
satisfies  for all $T>0$
\[ u_i \in L^2(Q_T), \; u_i \in L^\infty \big((0,T); H^{- 1}(\Omega) \big), \;
\int_0^T \int_\Omega u_1^\alpha u_2^\beta dx dt \leq C(T).
\]
2) 
For $\alpha= \beta =1$ and whatever are $\gamma, \; \delta$, System (\ref{eq:alphabetagammadelta}) has a non-negative weak solution which 
satisfies  for all $T>0$
\[u_i \in L^2(Q_T), \; u_i \in L^\infty \big((0,T); H^{- 1}(\Omega) \big), \;
\int_0^T \int_\Omega u_3^\gamma u_4^\delta dx dt \leq C(T).
\]
\end{theorem} 
\proof
We proceed in the same way as the previous theorem. Using the conservation laws for $\gamma u_1 +\alpha u_3$ and $\delta u_2+\beta u_4$ and Lemma~\ref{lm:mpg}, we obtain  $u_i$ is bounded  in $L^2(Q_T)$ independently of $n$.\qed

\begin{remark}
\label{rmk-alphabetagammadelta}
As in Remark~\ref{rmk-alphabetagamma}, we observe that\\
$\bullet$ when $\gamma=\delta=1$,  the specific form $u_1^\alpha u_2^\beta$ does not play any role, and any function $f(u_1,u_2) \geq 0$ satisfying $f(0,u_2)= f(u_1,0)=0$ will work ; \\
$\bullet$ when $\alpha =\beta= 1$, the specific form $u_3^\gamma u_4^\delta$ does not play any role, and any function $g(u_3,u_4) \geq 0$ satisfying $g(0,u_4)= g(u_3,0)=0$ will work
\end{remark}
\begin{remark}
As we see, the hypothesis \textbf{\em (QG)} is not used in the proofs of these last two theorems. The a priori controls given by the conservation laws allow us to conclude without the assumption of quadratic growth.
\end{remark}

%%%%%%%%%%%%%%%%%%%%%%%%%%%%%%%%%%%%%%%%%%%%
\subsection{A classical quadratic model with mass control }
\label{sec:quadratic2}
%-------------------------------------------
%%%%%%%%%%%%%%%%%%%%%%%%%%%%%%%%%%%%%%%%%%%%
%As a first example of the usefulness of our estimate (\ref{eq:mpg2}), 
%we show here how one may prove global existence of weak solutions for quadratic nonlinearities. In fact, this case  is of interest due to its relevance in many applications such as chemical reactions or Lotka-Volterra type systems.
%%Lotka-Volterra systems of the type described in [Leung] (see also [FHM]) and given as follows: 
Next, we establish global existence in time to another class of reaction-diffusion systems, which includes, among others, Lotka-Volterra systems. Namely, we consider the following system: for $t\geq 0, \; x \in \Omega \subset \R^N$, $1\leq i\leq m$, 
\begin{equation} \left\{
\begin{array}{ll}
	\dis\f{\p}{\p t} u_i(t,x)-  d_i\Delta u_i = f_i(u_1,\cdots,u_m) ,
	\\[5pt]
	\dis\f{\p u_i}{\p \nu} (t,x) =0  \text{ on } (0,+\infty)\times\p \Omega,
	\\[5pt]
	u_i(t=0)=u_i^0 \geq 0 ,
\end{array} \right.
\label{eq:quadContr}
\end{equation}
where for all $1\leq i\leq m$,   $d_i>0$ and $f_i:\R^m\rightarrow \R$ is locally Lipschitz.
\begin{theorem}[Quadratic growth] \label{quadgrowContr}
 Assume the  nonlinearities $f_i$ satisfy  \textbf{(P)}, \textbf{(M')} and \textbf{(QG)}. Then,  for all $u^0=(u_1^0,\cdots, u_m^0) \in (L^1(\Omega)^+\cap H^{- 1}(\Omega))^m$, System (\ref{eq:quadContr}) has a non-negative  weak solution which satisfies for all $T>0$
$$
u_i \in L^2(Q_T), \qquad u_i \in L^\infty \big((0,T); H^{- 1}(\Omega) \big) .
$$
\end{theorem} 

\noindent\textbf{Comment:} As we shall see, the crux is that mass control prevents the blow-up which may otherwise occur, for example, when the right-hand side has quadratic growth.\\

\proof 
For the  simplicity of the presentation, we consider the truncated system as before, but we drop the indexation by $n$.

We first prove the mass control. Again multiplying each $i$-th equation by $a_i$ and adding the $m$ equations   gives, by using \textbf{(M')}:
\begin{equation}\label{somme4}
\f{\p}{\p t} \big[\sum_{i=1}^m a_iu_i \big]- \Delta \big[ \sum_{i=1}^ma_id_iu_i \big]\leq  C_0\left(\sum_{i=1}^m a_i+ \sum_{i=1}^m a_iu_i\right).
\end{equation}
Integrating (\ref{somme4}) on $\Omega$ to obtain, with $B= C_0\dis \sum_{i=1}^m a_i$, 
\begin{equation}\label{integration1}
\f{d}{d t} \int_\Omega \big[\sum_{i=1}^m a_iu_i(t) \big]\leq B |\Omega|+ C_0\int_\Omega\big [\sum_{i=1}^m a_iu_i(t) \big]
\end{equation}
Let $T>0$. By integrating the previous  Gronwall inequality (\ref{integration1}), we  obtain,  thanks to non-negativity, for all $t\in [0,T]$
$$
\sum_{1\leq i\leq m}\|a_iu_i^n(t)\|_{L^1(\Omega)}\leq  e^{C_0 T}\left[\sum_{1\leq i\leq m} \|a_iu_{i}^0\|_{L^1(\Omega)}+ B|\Omega|T\right].
$$

The next step is to prove the quadratic estimate. We define  $U=\dis e^{-C_0 t} \sum_{i=1}^m a_iu_i$, $V= e^{-C_0 t} \dis\sum_{i=1}^m a_id_iu_i$. From equation~\eqref{somme4}, we compute 
\[
\f{\p}{\p t} U- \Delta V =B-F, \quad \text{with}\quad   F\geq 0.
\]
Therefore, we can apply Lemma~\ref {lm:mpg} and conclude a quadratic estimate (notice that $\int_0^T \langle F\rangle dt$ is controlled thanks to the mass control above) 
\[
\int_0^T \int_\Omega  UV \leq K(T) + \| U^0\|_{H^{-1}}
\]
from which we immediately deduce that the  $u_i$'s are bounded  in $L^2(Q_T)$.  Thus, the $f_i(u)$'s are bounded  in $L^1(Q_T)$ according the assumption~\textbf{(QG)}. Therefore, we can apply Theorem \ref{Lun}.
\qed

%%%%%%%%%%%%%%%%%%%%%%%%%%%%%%%%%%%%%%%%%%%%
\section{System with Porous Media diffusion }
\label{sec:porous}
%-------------------------------------------%

We continue our gallery of applications of Lemma~\ref{lm:mpg} with an example based on the porous medium equation. Because most of the literature about porous media uses Dirichlet boundary condition, we also do so. Consider the following system: for $t\geq 0, \; x \in \Omega \subset \R^N$, $1\leq i\leq m,$
\begin{equation} \label{snld}
\left\{
\begin{array}{ll}
	\dis\f{\p}{\p t} u_i(t,x)-  d_i\Delta u_i^{m_i} = f_i(u_1,\cdots,u_m) ,
	\\[10pt]
%	\dis\f{\p u_i}{\p \nu} (t,x) =0  \text{ on } (0,+\infty)\times\p \Omega,
    u_i(t,x) = 0 \text{ on } (0,+\infty)\times\p \Omega,
	\\[10pt]
	u_i(t=0)=u_i^0 \geq 0 ,
\end{array} \right.
\end{equation}
where  $d_i>0$, $m_i>1$ and $f_i:\R^m\rightarrow \R$ is locally Lipschitz-continuous.\\
%
%\vskip1mm
%\textcolor{blue}{Notice that the boundary conditions of this system are of Dirichlet type, not Neumann as previously in this article. We have not been able to carry over our proof to the Neumann case. However, existence for the system~\eqref{snld} was previously known with  only  $L^2$ initial data \cite{LaamPierAHP}[Corollary 2.8].}\\
%
\vskip1mm
Moreover, we assume that  the nonlinearities $f_i$ satisfy:\\
%
%$(H_1$) \textbf{quasi-positivity:}  $\forall 1\leq i\leq m$, 
%$$\forall r=(r_1,\cdots,r_m) \in [0,+\infty)^m,\; f_i(r_1,\cdots,r_{i-1},0,r_{i+1},\cdots, r_m)\geq 0 \;;$$
%%
%$(H_2)$ \textbf{mass dissipation:} there exists $(a_1,\cdots,a_n)\in (0,+\infty)^m$ such that 
%$$\forall r=(r_1,\cdots,r_m) \in [0,+\infty)^m,\quad\dis\sum_{i=1}^m a_if_i(r)\leq 0\; ;$$
%
\textbf{(SQG)}:  there exists $C>0$ and $\varepsilon>0$ such that $\forall 1\leq i\leq m$
$$ \forall \textbf{r}=(r_1,\cdots,r_m) \in [0,+\infty)^m,\quad|f_i(\textbf{r})|\leq C\big(1+\sum_{j=1}^m r_j^{m_j+1-\varepsilon}\big).$$
\begin{remark}
As is mentioned in \cite{LaamPierAHP}, the main point of the ``$-\varepsilon$'' in the above assumption is that it makes the nonlinearities $f_i^n(u^n)$ not only bounded in $L^1(Q_T)$, but \textit{uniformly integrable}. This is the main tool to pass to limit in the reactive terms.
\end{remark}
\noindent We now define what we mean by solution to our system (\ref{snld}).
\begin{definition}\label{veryweaksol}
Given $u_{i,0}\in L^1(\Omega)^+$ for $i=1,\cdots,m$.
We say that $u=(u_1,\cdots,u_m): (0,+\infty)\times\Omega \rightarrow (0,+\infty)^m$ is a global weak solution to system (\ref{snld}) if for all $i=1,\cdots,m$
\begin{equation}\label{weaksol}
\left\{
\begin{array}{l}
\text{\em for all } T>0,\\
u_i\in \mathcal{C}([0,T);L^1(\Omega)),\;u_i^{m_i}\in L^1(0,T;W_0^{1,1}(\Omega)), \; f_i(u) \in L^1(Q_T) \;\text{\em and }\;\forall \psi\in {\cal C}_T,\\
-\dis\int_\Omega \psi(0)u_i^0-\dis\int_{Q_T}\big[u_i\partial_t\psi+u_i^{m_i}\Delta\psi\big]=\dis\int_{Q_T}\psi\,f_i,
\end{array}
\right.
\end{equation}
where 
\begin{equation}\label{setofpsi} %
{\cal C}_T: =\{\psi:[0,T]\times\overline{\Omega}\rightarrow\R\;;\; \psi,\;\partial_t\psi,\;\partial^2_{x_ix_j}\psi \;{\rm are\; continuous},\; \psi=0\;{\rm on}\;\Sigma_T \;\text{\em and } \psi(T)=0\}.
\end{equation}
\end{definition}
The definition \ref{veryweaksol} corresponds to the notion of {\em very weak solution} in Definition 6.2 of \cite{Vaz}.\\ 
\begin{theorem}\label{quadgrowPM} %
 Assume that  the nonlinearities $f_i$ satisfy \textbf{(P)}, \textbf{(M)} and \textbf{(SQG)}. Then, for all $u^0=(u^0_1,\cdots, u^0_m) \in (L^1(\Omega)^+\cap H^{- 1}(\Omega))^m$, System (\ref{snld}) has a non-negative  weak solution which satisfies for all $T>0$
$$
u_i \in L^{m_i+1}(Q_T), \qquad u_i \in L^\infty \big((0,T); H_0^{- 1}(\Omega) \big).
$$
\end{theorem} 

\noindent \textbf{ Preliminary remark about {Theorem} \ref{quadgrowPM} and its proof:} Let us emphasize that, under the assumptions \textbf{(P)}, \textbf{(M)} and \textbf{ (SQG)}, the  existence of  global weak solutions to System~(\ref{snld}) subject to homogeneous  Dirichlet boundary conditions with initial data in $L^2(\Omega)$  was recently  established  by the first author and Pierre, see  \cite[Corollary 2.8]{LaamPierAHP}). The proof is based on a dimension-independent  $L^{m_i+1}(Q_T)$-estimate which heavily depends  on the $L^2$-norm of the initial data (see \cite[Theorem 2.7]{LaamPierAHP}), and then it is not valid for initial data in $L^1(\Omega)$. But, we are able to establish the same uniform bound in $L^{m_i+1}(Q_T)$ (see (\ref{bounded-in-L{mi+1}}))  with \textbf{\em initial data only in  $L^1(\Omega)^+\cap H^{-1}(\Omega)$} thanks to  the estimate (\ref{eq:mpg2}) of Lemma \ref{lm:mpg}. This is a \textit{ new result } in this  nonlinear degenerate setting and interesting for itself. 
Nevertheless, the hypothesis~\textbf{(SQG)} %``super growth" 
is crucial to pass to the limit.
\\  
Once the uniform estimate is  established,  the  demonstration by Laamri-Pierre in \cite{LaamPierAHP} applies {\em per se}. But,  for the sake of  completeness  and  for  the reader's convenience, we will give below  the main steps of the proof.\\

\proof We divide the proof in three steps.\\
$\bullet$ \textit{ First step:}  Let us  consider the following approximation of System (\ref{snld}) 
\begin{equation}\label{Systemn}
\left\lbrace\begin{array}{lll}
i=1,...,m,\\
%{\rm for\; all\;} T>0, u_i^n\in L^\infty(Q_T)\;\;u_i\geq 0,\;\; \varphi_i(u_i^n)\in %L^2(0,T;H_0^1(\Omega)),\\

\frac{\p u_i^n}{\p t}-d_i\Delta ((u_i^n)^{m_i})=f_i^n(u^n) &\text{ in } & Q=(0,+\infty)\times\Omega, 
\\[7pt]
 u_i^n =  0 &\text{ on }& \Sigma=(0,+\infty)\times\p \Omega,
\\[7pt]
u_i^n(0,.)=u_{i,0}^n\geq 0 & \text{ in } & \Omega ,\\
\end{array}
\right.
\end{equation}
where $u_{i,0}^n:= \inf\{u_{i,0},\; n\}$ and  $f_i^n:=\dis\frac{f_i}{1+\frac{1}{n}\sum_{1\leq j\leq m}|f_j|}$.\\
 For $i=1,\cdots,\,m$, $u_{i,0}^n\in L^\infty(\Omega)$ and converges to $u_{i,0}$ in $L^1(\Omega)$, $f_i^n$ is locally Lipchitz continuous and satisfies \textbf{(P), (M)} and \textbf{(SQG)}. Moreover
$\|f_i^n\|_{L^\infty(\Omega)}\leq n$. 
Therefore,  the approximate system~(\ref{Systemn}) has a non-negative bounded global solution  (see e.g  \cite[Lemma 2.3]{LaamPierAHP}). Moreover, this solution is regular in the sense  that
$$\text{for all } T>0,\; u_i^n \in L^\infty (Q_T) \cap \mathcal{C} ([0,T] ; L^1(\Omega)), \;  (u_i^n)^{m_i} \in L^2(0,T;H^1(\Omega))$$
and the equations  in (\ref{Systemn}) are satified a.e.\\
 Consequently, we have for $i=1,\cdots,\;m$
\begin{equation}\label{weakformuation}
\forall \psi\in {\cal C}_T,\\
-\dis\int_\Omega \psi(0)u_{i,0}^n-\dis\int_{Q_T}\big[u_i^n\partial_t\psi+(u_i^n)^{m_i}\Delta\psi\big]=\dis\int_{Q_T}\psi\,f_i^n.
\end{equation}
\textbf{Our goal is  to pass to the limit as $n\rightarrow +\infty$ in (\ref{weakformuation})}.\\
% We will do it in two steps.\\
$\bullet$ \textit{Second step }: 
 %This is where the main difficulty is concentrated. Actually, our method here, similar  to the one in   \cite[Corollary 2.8]{LaamPierAHP}, 
 First we will     establish a priori estimates and then we prove that  the $f_i^n(u^n)$ are uniformly bounded in $L^1(Q_T)$. For this, we multiply each $i$-th equation by $a_i$ and adding the $m$ equations to obtain 
\begin{equation}\label{somme}
\f{\p}{\p t} \big[\sum_{i=1}^ma_iu_i^n  \big]- \Delta  \big[\sum_{i=1}^md_ia_i(u_i^n)^{m_i}  \big]=\sum_{i=1}^m a_if_i^n(u^n).
\end{equation}
Thanks to  assumption  \textbf{(M)} and the definition of $f_i^n$, we obtain
 \begin{equation}\label{formeagreable}
 \f{\p}{\p t} U_n- \Delta V_n\leq 0 \text{ in } Q_T,\quad  V_n=0 \text{ on } \Sigma_T.
 \end{equation} 
where we  set  $U_n:=\dis\sum_{i=1}^ma_iu_i^n$ and $V_n: =\dis\sum_{i=1}^md_ia_i(u_i^n)^{m_i}$.
In order to apply Lemma~\ref{lm:mpg},  we first prove that  $<V_n(t)>$ is bounded.\\
To this end,  introduce $\theta$ solution of
\begin{equation}
-\Delta \theta= 1 \text{ in } \Omega,\quad \theta =0 \text{ in } \p\Omega.
\end{equation}
Now, integrating (\ref{formeagreable}) in time  leads to
\begin{equation}\label{integration-temps}
U_n(t,x) -\Delta \int_0^tV_n(s,x) ds\leq U_n(0,x)\leq U(0,x).
\end{equation}
Multilpying  (\ref{integration-temps}) by $\theta$ and integrating  over $\Omega$, to obtain
\begin{equation}
\int_\Omega \theta U_n(t,x)dx +\int_\Omega\theta\left(-\Delta[\int_0^t V_n(s,x)ds]\right
)\leq \int_\Omega\theta U(0,x) dx\leq \|\theta\|_{L^{\infty}(\Omega)}\|U_0\|_{L^{1}(\Omega)}
\end{equation}
After integration by parts, we obtain
\begin{equation} \label{fundamental}
\int_\Omega\int_0^t V_n(s,x) dsdx = \int_\Omega\int_0^t V_n(s,x)(-\Delta \theta) dsdx \leq \|\theta\|_{L^{\infty}(\Omega)}\|U_0\|_{L^{1}(\Omega)}
\end{equation}
which gives the desired estimate.\\
\noindent Thus,  according to  Lemma \ref{lm:mpg}, there exists  $C(T)>0$ independent of $n$ such that 
$$
\int_0^T \int_\Omega \big[\sum_{i=1}^m a_iu_i^n\big]\big[\sum_{i=1}^ma_id_i(u_i^n)^{m_i}\big]dtdx\leq C(T).
$$
Thanks to the non-negativity,  there exists $C>0$ independent of $n$ such that for all $i=1,\cdots\, m$, 
 \begin{equation}\label{bounded-in-L{mi+1}}
 \|u_i^n\|_{L^{m_i+1}(Q_T)}\leq C.
 \end{equation}
  Together with the assumption  \textbf{(SQG)}, the $f_i^n(u^n)$ are uniformly bounded in $L^1(Q_T)$, more precisely, there exists  a constant $C>0$ independent of $n$ such that
\begin{equation}\label{bound-of-fn}
\|f_i^n(u^n)\|_{L^1(Q_T)}\leq C. 
\end{equation}
In fact, $\{f_i^n(u^n)\}$ is not only bounded in $L^1(Q_T)$ but is also uniformly integrable on $Q_T$. Indeed, for all measurable set $E\subset Q_T$ with Lebesgue measure denoted by $|E|$, we have (recall that $|f_i^n|\leq|f_i|$)
$$\int_E|f_i^n(u^n)|\leq C\left[|E|+\int_E(u_i^n)^{m_i+1-\varepsilon}\right]\leq C\left[|E|+\left(\int_{Q_T}(u_i^n)^{m_i+1}\right)^{\frac{m_i+1-\varepsilon}{m_i+1}}|E|^{\frac{\varepsilon}{m_i+1}}\right].$$
Thanks to (\ref{bounded-in-L{mi+1}}),   $\int_E\sum_i|f_i^n(u^n)|$ may be made uniformly small by taking $|E|$ small enough. This is exactly the uniform integrability of the $f_i^n(u^n)$.\\ 
%
%Moreover, $f_i^n(u^n)$ converges a.e. to $f_i(u)$. Therefore, at least up to a subsequence, by Vitali's Lemma \ref{vitali}, we may deduce that $f_i^n(u^n)$ converges in $L^1(Q_T)$ for all $T<+\infty$ to $f_i(u_i)$. This implies that $u_i^n=S_{\varphi_i}\left(u_{i0}^n,f_i^n(u^n)\right)$ converges to $u_i=S_{\varphi_i}(u_{i0},f_i(u))$. 
%%%
%Finally, by the estimate (\ref{lukkari2}) in Lemma \ref{CompaciteLukkari} of $\nabla\varphi_i(u_i^n)$ in $L^\beta(Q_T)$ with $\beta>1$, it follows that $\varphi_i(u_i)$ is (at least) in $L^1(0,T;W^{1,1}_0(\Omega))$. This ends the proof of Corollary \ref{subquadratic}.
%$\blacksquare$\\

%%%%%%%%%%%%%%%%%%%%%%%%%%%%%%
\noindent$\bullet$ \textit{Third step:} Passage to the limit  as $n\rightarrow +\infty$ in (\ref{weakformuation})\\
%Now we are ready to pass to the limit in as $n\rightarrow +\infty$ in (\ref{weakformuation}).\\
%Since the $f_i^n$ are uniformly bounded in $L^1(Q_T)$, 
First, we apply the following compactness lemma wich allows us to extract  a converging subsequence from the $u_i^n$.
\begin{lemma}[Baras, \cite{Bar}]\label{baras}
Let $(w_0, H)\in L^1(\Omega)\times L^1(Q_T)$ and let $w$ be the solution of
$$w_t-d\Delta (w^q)=H \text{ \em in } Q_T,\quad  w=0 \text{ \em in } \Sigma_T, \quad w(0, .) =w_0.$$
The mapping $(w_0, H)\mapsto w$ is compact from $L^1(\Omega)\times L^1(Q_T)$ to $L^1(Q_T)$ for all $q>\dis\frac{(N-2)^+}{N}$.
\end{lemma}
Since the $f_i^n(u^n)$ is uniformly bounded in $L^1(Q_T)$, according to Lemma \ref{baras}, 
$\{u_i^n\}$ is relatively compact in  $L^1(Q_T)$.  Therefore, up to a subsequence, $\{u_i^n\}$ converges in  $L^1(Q_T)$ and a.e to some limit $u_i\in L^1(Q_T)$. Moreover, $f_i^n(u^n)$ is uniformly integrable   and converges a.e $f_i(u)$. Then, Vitali's theorem  implies that  $f_i(u) \in L^1(Q_T)$ and  $f_i^n(u^n)$ converges in $L^1(Q_T)$ to  $f_i(u)$.\\
Second, thanks to the $L^{m_i+1}$-estimate (\ref{bounded-in-L{mi+1}}), there exists a subsequence which converges  to $u_i$ in  $L^p(Q_T)$ for all $p<m_i+1$.\\
Now we can pass to the limit  as $n\rightarrow +\infty$ in the weak formulation~(\ref{weakformuation}) to obtain 
\begin{equation}
-\dis\int_\Omega \psi(0)u_i^0-\dis\int_{Q_T}\big[u_i\partial_t\psi+u_i^{m_i}\Delta\psi\big]=\dis\int_{Q_T}\psi\,f_i,
\end{equation}
\rm for all  $\psi \in {\cal C}_T$ , recall that 
\begin{equation*}\label{setofpsi2} %
 {\cal C}_T: =\{\psi:[0,T]\times\overline{\Omega}\rightarrow\R\;;\; \psi,\;\partial_t\psi,\;\partial^2_{x_ix_j}\psi \;{\rm are\; continuous},\; \psi=0\;{\rm on}\;\Sigma_T \;\text{\em and } \psi(T)=0\}.
\end{equation*}
To prove that $u=(u_1,\cdots,u_n)$ is solution to System (\ref{snld}) in the sense of definition \ref{veryweaksol}, it remains to establish that  for all  $T>0$,
$u_i\in \mathcal{C}([0,T);L^1(\Omega))$ and $ u_i^{m_i}\in L^1(0,T;W^{1,1}(\Omega))$ ; in fact, we will prove that $ u_i^{m_i}\in L^\beta(0,T;W_0^{1,\beta}(\Omega))$ for all $\beta\in [1,1+\dis\frac{1}{1+m_iN})$.
\begin{lemma}[Lukkari, \cite{Luk}]\label{CompaciteLukkari} % LEMME 3 ; 3.3 LUKKARI
Let $(w_0, H)\in L^1(\Omega)\times L^1(Q_T)$ and let $w$ be the solution of
\begin{equation}\label{eq:pm}
w_t-d\Delta (w^q)=H \text{ \em in } Q_T,\quad  w=0 \text{ \em in } \Sigma_T, \quad w(0, .) =w_0.
\end{equation}
Then, there exists $C>0$
\begin{equation}\label{lukkari1}
\int_{Q_T}|w|^{q\alpha}\leq C \;{\rm \;for \; \;} 0<\alpha< 1+\frac{2}{qN}, 
\end{equation}
\begin{equation}\label{lukkari2}
\int_{Q_T}|\nabla w^q|^\beta\leq C \;\;{\rm for }\;\;1\leq\beta< 1+\frac{1}{1+qN},
\end{equation}
\begin{equation}\label{lukkari3}
\left\|\dis\frac{\partial w}{\partial t}\right\|_{L^1(0,T; W^{-1,1}(\Omega)}\leq C.
\end{equation}
where $C=C\left(T,\alpha,\, \beta,\, q,\, \|w_0\|_{L^1(\Omega)},\, \|H\|_{L^1(Q_T)}\right)$.
\end{lemma}

For a proof, see Lukkari  \cite[Lemma 4.7]{Luk}.
%for the case $q> 1$ and Lukkari \cite[Lemma 3.5]{L5} for the case $\dis\frac{(N-2)^+}{N}<q< 1$. 
However, in this  reference, the proof is given with zero initial data, but with right-hand side a bounded measure. We may use the measure $\delta_{t=0}\otimes w_0\,dx$ to include the case of initial data $w_0$. We may also use the results in~\cite [Theorem  2.9]{APW}).  
\hfill \qed
\vskip1mm
Now let's come back to the proof of theorem.\\
Let $i\in \{1,\cdots, m\}$. According to the estimate (\ref{bound-of-fn}), we apply  the estimate (\ref{lukkari2}) to $i$-th equation of (\ref{snld}) to obtain 
% of the gradient in lemma \ref{CompaciteLukkari} 
 $\nabla (u_i^n)^{m_i}$ is bounded in the space $L^\beta(0,T;W_0^{1,\beta}(\Omega))$ for all $\beta\in [1,1+\dis\frac{1}{1+m_iN})$. These spaces being reflexive (for $\beta>1$), it follows that $\nabla (u_i)^{m_i}$ also belongs to these same spaces.\\
 According again (\ref{bound-of-fn}), we deduce from the estimate (\ref{lukkari2}) that $\dis\frac{\partial u_i}{\partial t}\in L^1(0,T; W^{-1,1}(\Omega))$  which implies in particular that 
 $u_i\in \mathcal{C}([0,T];L^1(\Omega))$.\qed

\begin{remark} % Remark 22
The paper by Gess, Sauer and Tadmor, \cite{GST19}, provides another route for the  compactness argument which can be applied directly to the $u_i^m$. It states, among other results, that solution of equation~\eqref{eq:pm} in the full space, with data in $L^1$ as in Lemma~\ref{CompaciteLukkari}, satisfy
\[
w \in L^{q}\big((0;T);W^{s,q}(\R^N) \big), \qquad s= \frac 2 q.
\]
and a similar space-time regularity result in fractional Sobolev spaces.
 \end{remark}

%%%%%%%%%
\appendix
%\section{Appendix}
%-----------------------------
%%%%%%%%%%%%

%-----------------------------------------------------------------
\section{First key  estimate with $L^1$ initial data}
\label{ap:KeyProof}
%------------------------------------------------------------------
%%%%%%%%%%%%%%%%%%%%%%%%%%%%%%%%%%%%%%%%%%%%
We prove the key Lemma~\ref{lm:mpg} and, for the sake of completeness, we recall the statement. We define
\[
 \langle F\rangle  = \f 1 {|\Omega |} \int_\Omega F(x) dx
\]
and, the $H^{-1}(\Omega)$ norm of $F$ as 
\[
\| F \|_{H^{-1}} = \| \nabla W \|_{L^2}, 
\] 
where $W$ solves
\[  \left\{
\begin{array}{ll}
\Delta W =F- \langle F\rangle , \quad x \in \Omega, 
\\[10pt]
	\f{\p W }{\p \nu} =0  \text{ on } \p \Omega.
\end{array} \right.
\]
%-------------------------------------
\begin{lemma}  [First key  estimate with $L^1$ data] \label{lm:mpgApp} 
Consider smooth functions $F, \; U: [0,+\infty)\times \Omega \rightarrow \R^+$ and $V: (0,+\infty)\times \Omega \rightarrow \R$ such that  $\dis\int_\Omega \Delta_x V(t,x)dx=0$ and $B$ a non-negative constant. Assume that  $U^0= U(t=0)\in L^1(\Omega)\cap H^{- 1}(\Omega)$  and that the differential relation holds
\begin{equation} \left\{
\begin{array}{ll}
	\f{\p}{\p t} U(t,x)-  \Delta_x V = B-F(t,x) \leq 0, \quad t\geq 0, \; x \in \Omega \subset \R^N, 
	\\[10pt]
	\f{\p U}{\p \nu}(t,x) =0  \;\text{\em  on } (0,+\infty)\times \p \Omega.
\end{array} \right.
\label{eq:mpg1App}
\end{equation}
 Then,  for some constant $C$ depending on $\Omega$, the inequality holds
\begin{equation}\label{eq:mpg2App}
\frac 1 2 \| U(T) \|_{H^{-1}}^2 + \int_0^T \int_\Omega  UV   \leq K(T) + \f 12 \| U^0\|_{H^{- 1}}^2 \quad \text{where}\quad   K=  \int_0^T \big[ C \langle F(t) \rangle + \langle V(t) \rangle \big]\; \big[B |\Omega| t+  \int_\Omega U^0 \big] dt .
\end{equation}
 \end{lemma}

Notice that,  with our assumptions,  $F$ is controlled in $L^1(Q_T)$ because
\[
 0 \leq  \langle U(T) \rangle  =  \langle U^0 \rangle  +B T  - \int_0^T  \langle F(s) \rangle  ds .
\]

\vskip2mm

\proof 
The proof is reminiscent from the lifting method introduced  in \cite{BloweyElliott91}. We write, substracting its average to equation~\eqref{eq:mpg1}, 
\[
\p_t (U- \langle U \rangle ) - \Delta V =  \langle F \rangle  - F .
\] 
Then, we solve with Neumann Boundary Condition
\[
- \Delta W  = U- \langle U \rangle ,  \qquad  \langle W \rangle =0  .
\]
Therefore, we find
\[
-\Delta \big[ \p_t W + V\big] \leq   \langle F \rangle  \qquad \text{(a constant)} ,
\]
and from elliptic theory, see \cite{ADN1, ADN2}, we conclude that
\beq
\p_t W+ V \leq C  \langle F \rangle  +  \langle V \rangle  .
\label{eq:fundamental_elliptic}
\eeq
Finally, multiplying by $U$,
\[
\int_\Omega \underbrace{\big[-\Delta W + \langle U \rangle \big]}_{=U}  \p_t W +\int_\Omega UV \leq \big[C  \langle F \rangle  +  \langle V \rangle \big] \int_\Omega U
\]
and, because $ \int_\Omega \p_t W =0$, this gives
\[
\f 12 \f{d}{dt} \int_\Omega |\nabla W (t)|^2 +\int_\Omega  UV(t) \leq \big[C  \langle F \rangle  +  \langle V \rangle \big] \big[ B |\Omega| t +\int_\Omega U^0\big] .
\]
The result follows after time integration.
\qed

\vskip2mm
\noindent \textbf{Comment :} An immediate variant of Lemma \ref{lm:mpgApp} holds true with the Dirichlet boundary condition and its proof is even simpler because subtracting the averages is useless and  the upper bound \eqref {eq:fundamental_elliptic} can be found thanks to explicit super-solutions.

%%%%%%%%%%%%%%%%%%%%%%%%%%%%%%%%%%%%%%%%%%%%
\section{A key estimate with $L^1$ source, no sign condition}
\label{sec:estimateA}
%-------------------------------------------
%%%%%%%%%%%%%%%%%%%%%%%%%%%%%%%%%%%%%%%%%%%%

 Lemma \ref{lm:mpgApp}  can be adapted to the case when the right-hand side $F$ does not have a sign. However, this loss of non-positivity is a major difficulty and this leads to the technical restriction $N\leq 3$. Because elliptic regularity is involved, the situation is more complicated and it is not clear whether the restriction $N\leq 3$ is  only technical or due to deeper phenomena. As we will see, it naturally appears in the proof.
\\

Again, for functions  $V, \; U, \; F: (0,+\infty)\times \Omega \rightarrow \R$ such that  $\dis\int_\Omega \Delta_x V=0$, we consider the  differential relation 

\begin{equation} \left\{
\begin{array}{ll}
	\f{\p}{\p t} U(t,x)-  \Delta V = F \in L^1, \quad t\geq 0, \; x \in \Omega \subset \R^N, 
	\\[5pt]
	\f{\p U}{\p \nu} (t,x) =0  \text{ on } (0,+\infty)\times \p \Omega,
	\\[5pt]
        U(t=0)=U^0\geq 0.
\end{array} \right.
\label{eq:mps1}
\end{equation}

%-------------------------------------
\begin{lemma} \label{lm:mps2} 
 Assume $N \leq 3$, $V \geq a \, U \geq 0$ for some $a>0$, $U^0\in L^1(\Omega)^+\cap H^{- 1}(\Omega)$, and  $\dis\int_\Omega \Delta_x V=0$. Then, it holds 
\begin{equation}
a \int_0^T \int_\Omega  U^2   \leq  C \left( \int_0^T [ \| F (t) \|_1 +\langle V(t) \rangle ] dt   \right)^2 + \| U^0\|_{H^{- 1}(\Omega)}^2.
\label{eq:mps2}
\end{equation}
\label{lm:mps} \end{lemma}

\proof 
As in Appendix~\ref{ap:KeyProof}, we first substract its average to equation~\eqref{eq:mps1} and find 
\[
\p_t (U- \langle U \rangle ) - \Delta V = F-  \langle F \rangle  .
\] 
Next, we solve with Neumann Boundary Condition
\[
- \Delta W  = U- \langle U \rangle , \qquad  \langle W \rangle =0  .
\]
Then, as in Appendix~\ref{ap:KeyProof},  we find 
\[
-\Delta \big[ \p_t W + V \big] = F-  \langle F \rangle   ,
\]
and from elliptic theory, \cite{ADN1, ADN2},  since $ 1 - \f 2 N < \f 1 2$ when $N \leq 3$,  
\[
\Phi:= \p_t W + V, \quad \text{satisfies}\quad   \| \Phi \|_{2} \leq C \|F\|_1 +  \langle V \rangle  .
\]
Therefore, multiplying again by $U$, 
\[
\int_\Omega\underbrace{[-\Delta W+ \langle U \rangle ]}_{=U}  \p_t W +\int  UV =  \int \Phi U \leq \| \Phi(t) \|_2   \| U(t) \|_2.
\]
Because $ \dis\int \p_t W =0$, this gives
\[
\f 12 \f{d}{dt} \int |\nabla W|^2 + a \int  U^2 \leq    \| \Phi(t) \|_2   \| U(t) \|_2.
\]
Therefore, we conclude after time integration.
\qed

%%%%%%%%%%%%%%%%%%%%%%%%%%%%%%%%%%%%%%%%%%%%
%-------------------------------------------
\section{Proof of Pierre's $L^2$ Lemma}
\label{Ap:Pierre}
%-------------------------------------------
%%%%%%%%%%%%%%%%%%%%%%%%%%%%%%%%%%%%%%%%%%%%

Integrating in time  the equation \eqref{eq:mpg1} with $V=AU\geq 0$, one finds, setting $W(t,x)= \int_0^t AU(s,x) ds$,
\[
U(t,x) - \Delta  W(t,x) \leq  U^0,
\]
We multiply by $AU= \f{\p W}{\p t}$ and find successively
\[
A U^2(t,x) -   \f{\p W}{\p t}\Delta  W(t,x) \leq A U U^0 ,
\]
\[
 \int_\Omega A U^2(t,x) +  \int_\Omega \f{\p \nabla W}{\p t}\nabla  W(t,x) \leq  \int_\Omega AU U^0 ,
\]
\[
 \int_\Omega A U^2(t,x) + \f 12 \f{d}{dt}  \int_\Omega |\nabla  W(t,x)|^2  \leq   \int_\Omega A U U^0 ,
\]
\[
\int_0^T \int_\Omega A U^2(t,x) + \f 12  \int_\Omega |\nabla  W(T,x)|^2  \leq  \left( \int_0^T \int_\Omega A (U^0)^2  \int_0^T \int_\Omega A U^2 \right)^{1/2}  .
\]
As a consequence, we arrive at the conclusion that 
\[
\int_0^T \int_\Omega A U^2(t,x)  \leq \sqrt {  \left\|   \int_0^T A \right\| _\infty }  \left( \int_\Omega (U^0)^2 \right)^{1/2} \left( \int_0^T \int_\Omega A U^2 \right)^{1/2}  .
\]
Comment: In this proof, the $L^2$ bound arises the time derivative, while in Lemma~\ref{lm:mpgApp} it stems from the Laplacian. 

%------------------------------------------------------------------------
\section{A priori bound for the porous medium equation}
\label{Ap:PM}
%------------------------------------------------------------------------

For the semi-linear porous medium equation, we need a consequence of Lemma~\ref{lm:mpgApp}  that we establish here. With the notations of section~\ref{sec:porous} and  for sake
of simplicity, we assume $a_1=\cdots=a_m=1$ and $d_1=\cdots=d_m=1$. By summing the $m$ equations, we obtain the problem 
\[
\f{\p}{\p t}  \sum_i u_i -\Delta \sum_i u_i^{m_i} = -F \in L^1(Q_T)  ,
\]
where $F:=-\dis\sum_i f_i \geq 0$.
\\
We prove that for some constant which depends on $T$, $\| u_i^0 \|_{L^1(\Omega)\cap H^{-1}(\Omega)}$ and  $\| F\|_{L^1(\Omega)}$, we have
\beq \label{pm:FundEstimate}
\int_0^T \int_\Omega   \sum_i u_i \; \sum_i u^{m_i} dx dt  \leq C \qquad \text{and } \; \int_0^T \int_\Omega  \sum_i u^{m_i} dx dt  \leq C .
\eeq
To prove these inequalities, we apply formula~\eqref{eq:mpg2} with $B=0$ and find
\[
\int_0^T \int_\Omega   \sum_i u_i \; \sum_i u^{m_i}    \leq C+ C \int_0^T   \int_\Omega   \sum_i u_i^{m_i} dx  dt .
\]
Therefore, using the H\"older's inequality
\[
\int_0^T \int_\Omega   \sup_i u_i \; \sum_i u^{m_i}    \leq C+ C  \sum_i  \left( \int_0^T \int_\Omega u_i^{m_i +1}  dx  dt\right)^{\f{m_i}{m_i +1} }
\]
and with $\al= \sup_i \f{m_i}{m_i +1}<1 $, 
\[ \bea
\dis \int_0^T \int_\Omega   \sup_i u_i \; \sum_i u^{m_i}  dx dt & \leq C+ C \dis  \sum_i  \left(\int_0^T \int_\Omega u_i^{m_i +1}  dx  dt\right)^{\al } 
\\
& \leq C+ C \dis  \left( \int_0^T \int_\Omega \sum_i   u_i^{m_i +1}  dx  dt\right)^{\al }
\\
&   \leq C+ C  \dis \left( \int_0^T \int_\Omega   \sup_i u_i \;  \sum_i   u_i^{m_i }  dx  dt\right)^{\al }
\eea
\]
This proves that $\int_0^T \int_\Omega   \sup_i u_i \; \sum_i u^{m_i}  dx dt$ is bounded and thus the result \eqref{pm:FundEstimate}.

%
%%%%%%%%%%%%%%%%%%%%%%%%%%%%%%%%%%%
%
%%%%%% BIBLIO %%%%%%%%%%%%%%%%%%%%%%
%
%%%%%%%%%%%%%%%%%%%%%%%%%%%%%%%%%%%%
%\pagestyle{myheadings}


\begin{thebibliography}{99}
\bibitem{APW}{B. Abdellaoui, I. Peral, M. Walias} : ``\textit{Some existence and regularity results for porous medium and fast equations with a gradient term}". Trans. Amer. Math. Soc. 367 (2015), no. 7, 4757--4791. 

%\bibitem{Amann}{H. Amann}, \textit{Dynamic theory of quasilinear parabolic systems.III. Global existence}, Math. Z., 202(1989), 219-250.
		
\bibitem{ADN1} S. Agmon, A.  Douglis, L.  Nirenberg :   ``\textit{Estimates near the boundary for solutions of elliptic partial
 differential equations satisfying general boundary conditions. I}".  Comm. Pure Appl. Math.  12 (1959), 623--727.

\bibitem{ADN2} S. Agmon, A.  Douglis, L.  Nirenberg :   ``\textit{Estimates near the boundary for solutions of elliptic partial
 differential equations satisfying general boundary conditions. II}". Comm. Pure Appl. Math.  17 (1964), 35--92.

\bibitem{Bar}{P. Baras} : ``\textit{Compacit\'{e} de l'op\'{e}rateur $f\mapsto u$ solution d'une \'{e}quation non lin\'{e}aire} $\frac{du}{dt} +Au \ni~f $". C.R.A.S. S\'{e}rie A, t. 286 (1978), 1113-1116.

\bibitem{BL} J. Bebernes and A. Lacey: ``\textit{Finite-time blow up for semilinear reactive-diffusive systems}". J. Diff. Equ., 95 (1992), 105-129.

%\bibitem{BR} S. Benachour and B. Rebiai, \textit{ Global classical solutions for reaction-diffusion systems with nonlinearities of exponential growth.} J. Evol. Equ., 10, Nr 3 (2010), 511-527.

\bibitem{BLMP} M. Bendahmane, Th. Lepoutre, A. Marrocco, B. Perthame : ``\textit{Conservative cross diffusions and pattern formation through relaxation}". J. Math. Pures et Appli. {\bf 92} 6 (2009), 651-667.

		
\bibitem{BloweyElliott91}  J. F. Blowey and C. M.   Elliott : ``\textit{The Cahn-Hilliard gradient theory for phase separation with nonsmooth free energy. I. Mathematical analysis}".  European J. Appl. Math.  2  (1991),  no. 3, 233--280.
%\bibitem{BoPi1} D.\ Bothe and M. Pierre, \textit{Quasi steady-state approximation for a reaction-diffusion system with fast intermediate.} J. Math. Anal. Appl.\ ., Vol. 368, Issue 1 (2010), 120-132.
 
%\bibitem{BoPi2} D.\ Bothe and M. Pierre, \textit{The instantaneous limit for reaction-diffusion systems with a fast irreversible reaction.},  Discrete Contin. Dyn.\ Syst.\ ., Serie S, 5 (2012), 49-59.

%\bibitem{BPR} D. Bothe, M. Pierre and G. Rolland, \textit{Cross-Diffusion Limit for a Reaction-Diffusion System with Fasr Reversible Reaction}, Comm. in PDE, {\bf 37} (2012), 1940-1966.


\bibitem{CanDesvFell} J.-A. Canizo, L. Desvillettes, K. Fellner : ``{\it Improved duality estimates and applications to reaction-diffusion equations}". Comm. Partial Differential Equations 39 (2014), 1185-204.
%
\bibitem{CapGoudVass} {C. Caputo, T. Goudon, A. Vasseur } : ``{\it Solutions of the 4-species quadratic reaction-diffusion system are bounded and $\mathcal{C}^\infty$-smooth, in any space dimension}".  Anal. PDE 12 (2019), no. 7, 1773--1804.

%
\bibitem{DFPV}  L. Desvillettes, K. Fellner, M. Pierre, J. Vovelle  :
``\textit{About Global Existence for Quadratic Systems of Reaction-Diffusion}".
Adv. Nonlinear Stud. {\bf 7} (2007), 491--511. 

%\bibitem{DiPL} R.J. Di Perna, P.-L. Lions: ``{\it On the Cauchy problem for Boltzmann equation: global existence and weak stability}". Ann. of Math., 130(2) (1989) 321-366.
  
\bibitem{FellLaam}{ K. Fellner, E.-H. Laamri } : ``\textit{Exponential decay towards equilibrium and global classical solutions for nonlinear reaction-diffusion systems}". J. Evol. Equ. \textbf{16} (2016), 681-704.

\bibitem{FellMorgTang}{ K. Fellner, J. Morgan, B. Q. Tang } : ``\textit{Global classical solutions to quadratic systems with mass control in arbitrary dimensions }". Preprint.



\bibitem{F} {W. Feng : ``\textit{Coupled system of reaction-diffusion equations and Applications in carrier facilitated diffusion}". Nonlinear Analysis, Theory, Methods and Applications \textbf{17}(3) (1991), 285--311.}

%\bibitem{fischer} J. Fischer, {\it Global existence of renormalized solutions to entropy-dissipating reaction-diffusion systems}, Arch. Ration. Mech. Anal., 218(1) (2015), 553-587.

%\bibitem {FHM} W.B. Fitzgibbon, S.L. Hollis, J.J. Morgan, \textit{Stability and Lyapunov Functions for Reaction-Diffusion Systems.} SIAM J. Math. Anal., Vol. 28, No3, pp. 595--610 (1997).

%\bibitem{FLM} W.E.\ Fitzgibbon, M.\ Langlais, J.J.\ Morgan, \textit{A degenerate reaction-diffusion system modeling atmospheric dispersion of pollutants.} J. Math. Anal. Appl. {\bf 307} (2005), 415--432.

\bibitem{GST19} B. Gess, J. Sauer, E. Tadmor : ``\textit{Optimal regularity in time and space for the porous medium equation}".  Preprint  	arXiv:1902.08632 (2019).


%\bibitem{GG} A. Glitzky and K. G\"{a}rtner, \textit{Energy estimates for continuous and discretized electro-reaction-diffusion systems}, Nonlinear Ana. Th. Methods Appl., Ser. A, 70, 2 (2009) 788-805.

\bibitem{GV10} T. Goudon and A. Vasseur : ``\textit{Regularity Analysis for Systems of Reaction-Diffusion Equations}", Annales Sci. ENS (4) 43 (2010), no. 1, 117--142.

%%\bibitem{HaWe} A.\ Haraux, F.B.\ Weissler, \textit{Non-uniqueness for a semilinear initial value problem.} Indiana Univ.\ Math.\ J.\ {\bf 31} (1982), 167--189. 

 %\bibitem{HY} A. Haraux and A. Youkana, \textit{On a result of K. Masuda concerning reaction-diffusion equations}, T\^ohoku Math. J. 40 (1988), 159-163 .



%%\bibitem{H} D. Henry, \underline{Geometric Theory of semilinear Parabolic Equations}, Lecture Notes in Math., Vol. 840, Springer-Verlag, 1981.

%%\bibitem{HLV} M.A. Herrero, A.A. Lacey and J.L Vel\`azquez, \textit{Global existence for reaction-diffusion systems modelling ignition}, Arch. Rat. Mech. Anal., 142 (1998), 219-251.



%\bibitem{HMP} S.L. Hollis, R.H. Martin and M. Pierre : \textit{Global existence and boundedness in reaction-diffusion systems}, SIAM J. Math. Ana. 18 (1987), 744-761.

%\bibitem{K} J.I. Kanel, \textit{Cauchy's problem for semilinear parabolic equations with balance law}, Differentsial'naya Uravneniya, 20 (1984), 1753-1760.

%\bibitem{K1}  J.I. Kanel, \textit{The solvability on the whole of reaction-diffusion systems with condition of balance}, Differentsial'naya Uravneniya, 26 (1990), 448-458.

%\bibitem{KK} J.I. Kanel and M. Kirane, \textit{Global Solutions of Reaction-Diffusion Systems with a Balance Law and Nonlinearities of Exponential Growth}, J. Diff. Equ., 165, (2000), 24-41.

%\bibitem{Kou} S. Kouachi, \textit{Existence of global solutions to reaction-diffusion systems with nonhomogeneous boundary conditins via a Lyapounov functional}, Electron. J. Dif. Eqns., 88 (2002), 1-13.

%\bibitem{KouYou} S. Kouachi and A. Youkana, \textit{Global existence for a class of reaction-diffusion systems.} Bull. Pol. Acad. Sci. Math., {\bf 49}, Nr3 (2001).



\bibitem{Laamthese} E.-H.  Laamri : ``\textit{Existence globale pour des syst\`emes de r\'eaction-diffusion dans $L^1$}". Ph.D thesis, universit\'e de Nancy 1, France, 1988.

\bibitem{LaamAAM} E.-H. Laamri : ``\textit{Global existence of classical solutions for a class of reaction-diffusion systems}". Acta Applicandae Mathematicae, {\bf 115} 2 (2011), 153--165.


\bibitem{LaamPierAHP}{E.-H.  Laamri,  M. Pierre } : ``\textit{ Global existence for reaction-diffusion systems with nonlinear diffusion and control of mass}". Ann. Inst. H. Poincaré Anal. Non Linéaire 34 (2017), no. 3, 571--591.

\bibitem{LaamPierM3AS}{E.-H Laamri,  M. Pierre }: ``\textit{Stationary reaction-diffusion systems in $L^1$}". Mathematical Models and Methods in Applied Sciences
Vol. 28, No. 11 (2018) 2161--2190.

\bibitem{LSU}{O. A. Ladyzenskaya, V. A. Solonnikov, N. N. Uralceva : \underline{``Linear and quasilinear equations } \underline{of parabolic type}". Transl. Math. Monographs, \textbf{23}, A.M.S., Providence, R.I. 1968.}



%\bibitem{Leung} A. Leung, \underline{Systems of Nonlinear Partial Differential Equations,} Kluwer Academic Publ. Boston, 1989.

%\bibitem{reference7} G.M. Lieberman, {\it Second Order Parabolic Differential Equations}, World Scientific, Singapore, 1996. 

\bibitem{Luk}{T. Lukkari} : ``\textit{The porous medium equation with measure data}". J. Evol. Equ., {\bf 10} (2010), 711-729.


%\bibitem{MPa} T. Manteuffel and S. Parter, \textit{Preconditioning and boundary conditions}, SIAM J. Numer. Anal. 27 (1990), 656-694.



 \bibitem{MP1} R. H. Martin, M. Pierre : ``\textit{Nonlinear reaction-diffusion systems}". In Nonlinear Equations in the Applied Sciences, W.F. Ames and C. Rogers ed., Math. Sci. Eng. 185, Ac.
Press, New York, 1991.

\bibitem{MP2} R. H. Martin,  M. Pierre : ``\textit{Influence of mixed boundary conditions in some
reaction-diffusion systems}". Proc. Roy. Soc. Edinburgh, section A 127 (1997), 1053-1066 .

%\bibitem{M} K. Masuda, \textit{On the global existence and asymptotic behavior of reaction-diffusion equations}, Hokkaido Math. J. 12 (1983), 360-370.

%\bibitem{MGR} J.S. McGough and K. L. Riley, \textit{A priori bounds for reaction-diffusion systems arising in chemical and biological dynamics.} Appl. Math. Comput.,  {\bf 163}, No1 (2005), 1--16.

\bibitem{Mo89} J. Morgan : ``\textit{Global existence for semilinear parabolic systems}". SIAM J. Math. Ana. 20 (1989), 1128-1144.

%\bibitem{Mo90} J.\ Morgan, \textit{Boundedness and decay results for reaction-diffusion systems.} SIAM J. Math. Anal. {\bf 21} (1990), 1172--1189.

%\bibitem{Mu} J.D. Murray, \underline{Mathematical Biology.} Biomath. texts, Springer, 1993. %

%\bibitem{Oth} H.G Othmer, F.R. Adler, M.A. Lewis and J. Dallon, eds, \underline{Case studies in Mathematical Modeling-Ecology, Physiology and Cell Biology.} Prentice Hall, New Jersey, 1997.

%\bibitem{Pao} C.V. Pao, \textit{Nonlinear Parabolic and Elliptic equations.} Plenum, New York, 1992.

\bibitem{BPerth} B.  Perthame : ``\underline{Parabolic equations in biology: 
growth, reaction, movement and diffusion}". Lecture Notes on Mathematical Modelling in the Life Sciences. Springer, Berlin, 2015.

%\bibitem{Pierre87} M. Pierre, {\em An $L^1$-method to prove global existence in some reaction-diffusion systems}, in "Contributions to Nonlinear Partial Differential Equations", Vol.II, Pitman Research notes, 155, J.I. Diaz and P.L. Lions ed., (1987) 220-231.

\bibitem{Pier2003} M. Pierre : ``\textit{Weak solutions and supersolutions in $L\sp 1$ for reaction-diffusion systems}". J. Evol. Equ.  \textbf{3},  (2003), no. 1, 153--168. 

\bibitem{Pier2010} M. Pierre : ``\textit{Global Existence in Reaction-Diffusion Systems with Control of Mass: a Survey}". Milan J. Math. Vol. 78 (2010), 417-455.

\bibitem{PierRoll} M. Pierre, G. Rolland : ``\textit{Global existence for a class of quadratic reaction-diffusion system with nonlinear diffusion and $L^1$ initial data}". J. Nonlinear Analysis TMA. Volume 138 (2017), 369--387.

\bibitem{PS1997} M. Pierre, D. Schmitt : ``\textit{Blow-up in Reaction-Diffusion Systems with
Dissipation of Mass}". SIAM J. Math. Ana. 28, No2 (1997), 259-269.

\bibitem{PS2000} M. Pierre,  D. Schmitt : ``\textit{Blow-up in 
reaction-diffusion systems with dissipation of mass}". 
SIAM Rev.  \textbf{42},  (2000),  pp. 93--106 (electronic). 

\bibitem{PierSuzYam2019} M. Pierre, T. Suzuki,  Y. Yamada : ``\textit{Dissipative reaction diffusion systems with quadratic growth}". Indiana Univ. Math. J. 68 (2019), no. 1, 291--322. 

\bibitem{Rot} F. Rothe : \underline{``Global Solutions of Reaction-Diffusion Systems}". Lecture Notes in Mathematics, Springer, Berlin,
(1984).

\bibitem{QS} {P. Quittner, Ph. Souplet : \underline{``Superlinear Parabolic Problems: Blow-up, Global Existence and Steady} \underline{States}". Advanced Texts, Birkh\"{a}user, (2007).    }


%\bibitem{S2}{D. Schmitt: \textit{Existence globale ou explosion  pour des syst\`emes de r\'{e}action-diffusion avec contr\^{o}le de masse}. Th\`ese. Univ. de Nancy 1, 1995.}

%%\bibitem{Si} Simon, J. (1987). \textit{Compact Sets in the Space $L^p(0,T;B)$}, Annali di Matematica pura ed applicata (IV){\bf 146}: 65--96.

\bibitem{Soup18} Ph. Souplet : ``\textit{Global existence for reaction-diffusion systems with dissipation of mass and quadratic growth}".
 J. Evol. Equ. 18 (2018), no. 4, 1713--1720.

%\bibitem{T} A.M. Turing, \textit{The chemical basis of morphogenesis.} Philos. Trans. R. Soc. London Ser. B {\bf 237} (1952), 37--72.

\bibitem{SuzYam}  T. Suzuki, Y. Yamada :  ``\textit{Global-in-time behavior of Lotka-Volterra system with diffusion: skew-symmetric case}". Indiana Univ. Math. J. 64 (2015), no. 1, 181--216. 

\bibitem{Vaz}{ J. L. Vazquez  : } \underline{``\textit{The porous medium equation Mathematical Theory}"}. Oxford Mathematical Monographs, 2006.


\end{thebibliography}
\end{document}